\documentclass[12pt]{amsart}
\usepackage{amsmath,amscd,amssymb,amsfonts,graphics}
\newcommand{\h}{\hbox}
\newcommand{\q}{\quad}

\newcommand{\bs}{\par\bigskip}
\newcommand{\ms}{\par\medskip}
\newcommand{\sk}{\par\smallskip}
\newcommand{\bsn}{\par\bigskip\noindent}
\newcommand{\msn}{\par\medskip\noindent}

\newcommand{\ges}{\geqslant}
\newcommand{\les}{\leqslant}
\newcommand{\1}{\hskip1pt}
\newcommand{\mcap}{\hbox{$\bigcap$}}
\newcommand{\mcup}{\hbox{$\bigcup$}}
\newcommand{\msum}{\hbox{$\sum$}}
\newcommand{\mopl}{\hbox{$\bigoplus$}}
\newcommand{\mprod}{\hbox{$\prod$}}

\newcommand{\D}{{\mathcal D}}
\newcommand{\Lc}{{\mathcal L}}
\newcommand{\OO}{{\mathcal O}}
\newcommand{\PP}{{\mathbb P}}
\newcommand{\Q}{{\mathbb Q}}
\newcommand{\C}{{\mathbb C}}
\newcommand{\N}{{\mathbb N}}
\newcommand{\R}{{\mathbb R}}
\newcommand{\RR}{{\mathbf R}}
\newcommand{\Z}{{\mathbb Z}}
\newcommand{\ob}{{\mathbf 1}}
\newcommand{\Gr}{{\rm Gr}}
\newcommand{\al}{\alpha}
\newcommand{\be}{\beta}
\newcommand{\ga}{\gamma}
\newcommand{\Ga}{\Gamma}
\newcommand{\De}{\Delta}
\newcommand{\la}{\lambda}
\newcommand{\si}{\sigma}
\newcommand{\ep}{\varepsilon}

\newcommand{\alt}{\widetilde{\alpha}}

\newcommand{\bt}{\widetilde{b}}

\newcommand{\V}{\widetilde{V}}
\newcommand{\Ht}{\widetilde{H}}
\newcommand{\M}{{}\,\widetilde{\!M}{}}
\newcommand{\Pt}{{}\,\widetilde{\!P}{}}
\newcommand{\Om}{\Omega}
\newcommand{\om}{\omega}
\newcommand{\dd}{\partial}
\newcommand{\ddd}{{\rm d}}
\newcommand{\ee}{{\bf e}}
\newcommand{\mm}{{\mathfrak m}}
\newcommand{\ib}{{\mathbf i}}
\newcommand{\Ff}{F_{\!f}}
\newcommand{\Gp}{\Gamma_{\!+}}
\newcommand{\LF}{{\rm LF}}
\newcommand{\LFp}{{\rm LF}_{\!+}}
\newcommand{\Spx}{{\rm Supp}_{(x)}}
\newcommand{\Sp}{{\rm Sp}}
\newcommand{\HI}{{\rm HI}}
\newcommand{\Tj}{{\rm Tj}}
\newcommand{\bl}{\bigl}
\newcommand{\br}{\bigr}

\newcommand{\ssb}{\raise.15ex\h{${\scriptscriptstyle\bullet}$}}
\newcommand{\ssc}{\,\raise.15ex\h{${\scriptstyle\circ}$}\,}
\newcommand{\onto}{\mathop{\rlap{$\to$}\hskip2pt\hbox{$\to$}}}
\newcommand{\into}{\hookrightarrow}
\newcommand{\simto}{\,\,\rlap{\hskip1.5mm\raise1.4mm\hbox{$\sim$}}\hbox{$\longrightarrow$}\,\,}
\begin{document}
\title[Hodge ideals and spectrum]
{Hodge ideals and spectrum of isolated\\hypersurface singularities}
\author[S.-J. Jung]{Seung-Jo Jung}
\address{S.-J. Jung : Department of Mathematics Education, Jeonbuk National University, Jeonju, 54896, Korea}
\email{seungjo@jbnu.ac.kr}
\author[I.-K. Kim]{In-Kyun Kim}
\address{I.-K. Kim : Department of mathematics, Yonsei University, 50 Yonsei-Ro, Seoul 03722, Korea}
\email{soulcraw@gmail.com}
\author[M. Saito]{Morihiko Saito}
\address{M. Saito : RIMS Kyoto University, Kyoto 606-8502 Japan}
\email{msaito@kurims.kyoto-u.ac.jp}
\author[Y. Yoon]{Youngho Yoon}
\address{Y. Yoon : Department of Mathematics, Chungnam National University, 99 Daehak-ro, Daejeon 34134, Korea}
\email{mathyyoon@gmail.com}
\thanks{This work was partially supported by BK21 PLUS SNU Mathematical Sciences Division and the National Research Foundation of Korea(NRF) grant funded by the Ministry of Science, ICT and Future Planning (the first author: NRF-2018R1D1A1B07046508, the second author: NRF-2020R1A2C4002510, and the fourth author: NRF-2020R1C1C1A01006782). The third author is partially supported by JSPS Kakenhi 15K04816.}
\begin{abstract} We introduce Hodge ideal spectrum for isolated hypersurface singularities to see the difference between the Hodge ideals and the microlocal $V$-filtration modulo the Jacobian ideal. Via the Tjurina subspectrum, we can compare the Hodge ideal spectrum with the Steenbrink spectrum which can be defined by the microlocal $V$-filtration. As a consequence of a formula of Mustata and Popa, these two spectra coincide in the weighted homogeneous case. We prove sufficient conditions for their coincidence and non-coincidence in some non-weighted-homogeneous cases where the defining function is semi-weighted-homogeneous or with non-degenerate Newton boundary in most cases. We also show that the convenience condition can be avoided in a formula of Zhang for the non-degenerate case, and present an example where the Hodge ideals are not weakly decreasing even modulo the Jacobian ideal.
\end{abstract}
\maketitle
\centerline{\bf Introduction}
\bsn
M.~Musta\c{t}\v{a} and M.~Popa (\cite{MP2}, \cite{MP3}) recently defined {\it Hodge ideals} $I_p(D)\subset\OO_X$ for $\Q$-divisors $D=\msum_k\,\al_kZ_k$ on smooth varieties $X$. Here the $Z_k$ are reduced and irreducible, and we assume $\al_k\in(0,1]$. These can be extended naturally to the analytic case, see (2.1) below. In this introduction, we assume that $Z:=\mcup_k\,Z_k$ is {\it irreducible,} that is, $Z=Z_1$, and put $\al=\al_1$.
\sk
Hodge ideals provide a quite intersecting refinement of {\it multiplier ideals}, especially in the case the minimal exponent $\alt_Z$ is at least 1, where the classical multiplier ideals become powers of the ideal of $Z\subset X$. It is shown there that $I_p(\al Z)$ coincides with the {\it microlocal $V$-filtration} $\V^{\al+p}\OO_X$ {\it modulo} $(f)$, where $f$ is a local defining function of $Z\subset X$, see also (2.4.7) below. (In the case $\al=1$, this was shown in \cite[Theorem 1]{hi}.) However, the relation between $I_p(\al Z)$ and $\V^{\al+p}\OO_X$ without modulo $(f)$ seems rather complicated in general, see the above papers of Musta\c{t}\v{a} and Popa as well as \cite{JKY}, \cite{hi}, \cite{Zh}, etc.
To see their difference without modulo $(f)$, it seems then interesting to compare these {\it modulo the Jacobian ideal} $(\dd f)\subset\OO_X$ generated by the partial derivatives $f_i:=\dd f/\dd x_i$, where $x_1,\dots,x_n$ are local coordinates of $X$ with $n=\dim X$.
\sk
From now on, we assume that $Z$ has an {\it isolated singularity\1} at $0\in X$ so that $\C\{x\}/(\dd f)$ is {\it finite-dimensional,} where $\C\{x\}=\C\{x_1,\dots,x_n\}=\OO_{X,0}$. In \cite{St2} (see also \cite{St3}), Steenbrink defined the {\it spectrum}
$$\Sp_f(t)=\msum_{i=1}^{\mu_f}\,t^{\al_{f,i}}\,\in\,\Z[t^{1/e}],$$
using the mixed Hodge structure on the vanishing cohomology together with the monodromy, where $\mu_f$ is the Milnor number, and $e$ is a positive integer related to the monodromy (more precisely, $T_s^e={\rm id}$ with $T=T_sT_u$ the Jordan decomposition of the monodromy), see also (1.1) below. The positive rational numbers $\al_{f,i}$ ($1\les i\les\mu_f$) are assumed weakly increasing, and are called the {\it exponents\1} or {\it spectral numbers\1} of $f$.
\sk
By \cite{SS}, \cite{Va1}, the spectrum $\Sp_f(t)$ can be defined also as the Hilbert-Poincar\'e series of the finite-dimensional filtered vector space
$$\Om^n_f:=\Om_{X,0}^n/\ddd f{\wedge}\Om_{X,0}^{n-1},$$
so that
$$\#\{i\mid\al_{f,i}=\be\}=\dim_{\C}\Gr_V^{\be}\Om_f^n\q(\be\in\Q).
\leqno(1)$$
where $V$ is the quotient filtration of the $V$-filtration on the Brieskorn lattice $H''_f$ in \cite{Br}, see (1.2) below.
\sk
Let $I_p(\al Z)\subset\C\{x\}$ be the Hodge ideals for $\al\in(0,1],\,p\in\N$, see (2.1) below. Since $Z$ has an isolated singularity at 0 and the $I_p(\al Z)$ are coherent, these are $\mm_{X,0}$-{\it primary ideals,} that is, $I_p(\al Z)\supset\mm_{X,0}^k$ for some $k\in\Z_{>0}$ (depending on $\al,p$) with $\mm_{X,0}\subset\C\{x\}$ the maximal ideal, and $\C\{x\}/I_p(\al Z)$ is finite-dimensional, see Remark~(2.1)\,(ii) below.
The {\it Hodge ideal spectrum\1} $\Sp_f^{\HI}(t)$ is defined as the Hilbert-Poincar\'e series of the finite-dimensional filtered vector space
$$(\Om_f^n,V_{\HI})\cong(\C\{x\}/(\dd f),V_{\HI}),$$
with $V_{\HI}$ defined by
$$V_{\HI}^{\be}\bl(\C\{x\}/(\dd f)\br):=\msum_{\al+p\ges\be}\,I_p(\al Z)\,\,\,\,\h{mod}\,\,\,\,(\dd f),
\leqno(2)$$
so that
$$\Sp_f^{\HI}(t)=\msum_{i=1}^{\mu_f}\,t^{\1\al^{\HI}_{f,i}}\q\h{with}\q \#\{i\mid\al^{\HI}_{f,i}=\be\}=\dim_{\C}\Gr_{V_{\HI}}^{\be}\bl(\C\{x\}/(\dd f)\br).$$
Here the $\al^{\HI}_{f,i}$ are assumed weakly increasing.
\sk
The above definition of Hodge ideal spectrum differs from the one in \cite{JKY}, where $V_{\HI}^{\be}\bl(\C\{x\}/(\dd f)\br)$ was defined by $I_p(\al Z)$ mod $(\dd f)$ for $\al{+}p=\be$ with $\al\in(0,1]$, $p\in\N$ {\it without\1} taking the above summation. However, the Hodge ideals $I_p(\al Z)$ mod $(\dd f)$ are not necessarily {\it weakly decreasing\1} for $\al\in(0,1]$ (with $p\ges 1$ fixed), see Example~(4.2) below. (Without taking mod $(\dd f)$, this was observed in \cite[Example 10.5]{MP2}, \cite[Example 4.6]{Zh}.)
Note that $I_0(\al Z)=\V^{\al}\C\{x\}$ for $\al\in(0,1]$, hence $\al_{f,i}=\al^{\HI}_{f,i}$ if $\al_{f,i}<1$, see \cite{MP2}, \cite{MP3}. In particular, we have $\Sp_f^{\HI}(t)\in\Z[t^{1/{e'}}]$, where $e'$ might be different from $e$ in general.
\sk
We can define also the {\it Tjurina subspectrum} $\Sp_f^{\Tj}(t)$ by
$$\Sp_f^{\Tj}(t)=\msum_{i=1}^{\tau_f}\,t^{\1\al^{\Tj}_{f,i}}\q\h{with}\q \#\{i\mid\al^{\Tj}_{f,i}=\be\}=\dim_{\C}\Gr_{V_{\HI}}^{\be}\bl(\C\{x\}/(\dd f,f)\br),$$
where $\tau_f$ is the Tjurina number of $f$, and the $\al^{\Tj}_{f,i}$ are assumed weakly increasing.
This gives a link between $\Sp_f(t)$ and $\Sp_f^{\HI}(t)$. Indeed, there is a subset $I\subset\{1,\dots,\mu_f\}$ such that $|I|=\tau_f$ and
$$\Sp_f^{\Tj}(t)=\msum_{i\in I}\,t^{\1\al^{\HI}_{f,i}},\q\Sp_f^{\HI}(t)-\Sp_f^{\Tj}(t)=\msum_{i\notin I}\,t^{\1\al^{\HI}_{f,i}},
\leqno(3)$$
since the canonical projection $\C\{x\}/(\dd f)\onto\,\C\{x\}/(\dd f,f)$ is strictly compatible with $V_{\HI}$.
(The geometric meaning of the Tjurina {\it subspectrum\1} does not seem to be clear.)
\sk
On the other hand, the $V$-filtration on $\Om_f^n\cong\C\{x\}/(\dd f)$ used in (1) coincides with the quotient filtration of the microlocal $V$-filtration $\V$ on $\C\{x\}$ under the canonical surjection
$$\C\{x\}\onto\,\C\{x\}/(\dd f),$$
see Proposition~(1.4) below. Then the coincidence of the Hodge ideals and the microlocal $V$-filtration modulo $(f)$ mentioned above implies that
$$\#\{i\mid\al^{\Tj}_{f,i}=\be\}=\dim_{\C}\Gr_{\V}^{\be}\bl(\C\{x\}/(\dd f,f)\br),$$
hence there is a subset $J\subset\{1,\dots,\mu_f\}$ such that $|J|=\tau_f$ and
$$\Sp_f^{\Tj}(t)=\msum_{i\in J}\,t^{\1\al_{f,i}},\q\Sp_f(t)-\Sp_f^{\Tj}(t)=\msum_{i\notin J}\,t^{\1\al_{f,i}}.
\leqno(4)$$
(So $\Sp_f^{\Tj}(t)$ is called a {\it subspectrum}.) In particular, we get the following (see also \cite{JKY}):
$$\Sp_f(t)=\Sp^{\HI}_f(t)=\Sp_f^{\Tj}(t)\,\,\,\h{if $f$ is weighted homogeneous.}
\leqno(5)$$
The last assumption is equivalent to that $\mu_f=\tau_f$, see \cite{SaK}. The partial converse of (5) with second equality forgotten does {\it not\1} necessarily hold as is seen by Theorem~1 below.
\sk
In our main theorems, we will often assume the following:
$$\aligned&\h{The function $f\in\C\{x\}$ is {\it semi-weighted-homogeneous\1} with weights $w_i$,}\\&\h{or $f$ has {\it non-degenerate\1} Newton boundary.}\endaligned
\leqno{\rm (A)}$$
Here we do not have to assume that $f$ is convenient, see (1.6--7) and Theorem~(2.7) below for more details. Note that $Z$ has an isolated singularity as is assumed above.
When we consider condition~(A), the coordinates $x_1,\dots,x_n$ are {\it fixed,} and cannot be replaced easily unless it is done with enough care.
\sk
In the main theorems we will sometimes assume that $f$ {\it is not a double point,} that is, $f\in\mm_{X,0}^3$. In this case, set
$$\ga_f:=\max\bl\{\be\in\Q\mid V^{\be}(\Om_f^n/\mm_{X,0}^2\Om_f^n)\ne 0\br\},$$
where $V$ is as in (1). This is compatible with the definitions of $\ga_f$ in the case condition~(A) is satisfied (see (1.6--7) below) assuming also $f\in\mm_{X,0}^3$, see \cite[Remark (ii) in Section 4.1]{mos}, \cite{exp}, \cite[Proposition 3.2]{mic}, \cite{VK} (and also (1.6--7) below).
\sk
If $f\in\mm_{X,0}^3$ or condition~(A) is satisfied, set
$$\ep_f:=\ga_f+1-\al_{f,\,\mu_f}=2\al_{f,1}-(n{-}1)+(\ga_f{-}\al_{f,1}).
\leqno(6)$$
\sk
The following theorem says that in some special cases this $\ep_f$ determines whether the Hodge ideal spectrum coincides with the Steenbrink spectrum (see (3.1) below).
\msn
{\bf Theorem~1.} {\it Assume condition~{\rm(A)} is satisfied, $f$ is not a double point, that is, $f\in\mm_{X,0}^3$, and moreover
$$f\Om_f^n=V^{\al_{f,\,\mu_f}}\Om_f^n\q(\h{in particular}\,\,\,\,\tau_f=\mu_f-1).
\leqno(7)$$
Then $\,\Sp_f^{\HI}(t)\ne\Sp_f(t)\,$ if and only if $\ep_f>0$. More precisely, we have}
$$\al_{f,\,\mu_f}^{\HI}-\al_{f,\,\mu_f}=\max\bl(\ep_f,0\br).
\leqno(8)$$
\ms
Condition (7) holds if $\tau_f=\mu_f-1$ and $\al_{f,\tau_f}-\al_{f,1}<1$ (although the last two conditions also imply that $\ep_f>0$).
As a corollary of Theorem~1, we see that the first equality of (5) can hold even in the {\it non-weighted-homogeneous} case if $\ep_f\les 0$. It is not difficult to construct examples with $\ep_f$ vanishing or positive or negative for $n=2$, see Example~(4.1) below. This also shows that the $\exp\bl(-2\pi i\1\al_{f,j}^{\HI}\br)$ ($j\in[1,\mu_f]$) can be different from the set of Milnor monodromy eigenvalues (answering a question of the referee), see Remark~(4.1) below.
\sk
As for the last equality in the definition of $\ep_f$, note that $\al_{f,\,\mu_f}=n-\al_{f,1}$ by the {\it symmetry\1} of spectral numbers, see \cite{St2} (and (1.1.2) below). This symmetry also implies that $2\al_{f,1}\les n$. In the case $\mu_f\ne\tau_f$, we have the inequality
$$2\al_{f,1}\les n{-}1,\q\h{or equivalently,}\q 2\al_{f,\,\mu_f}\ges n{+}1,
\leqno(9)$$
using (1.2.6) below. The condition $2\al_{f,1}>n{-}1$ is equivalent to that $Z$ has a singularity of type $A$, $D$, $E$ (at least for $n\les 3$) according to the classification theory of holomorphic functions with isolated singularities.
\sk
For the proof of Theorem~1, note that the assertion is equivalent to the equality
$$\al_{f,\,\mu_f}^{\HI}=\max\bl(\ga_f{+}1,\,\al_{f,\,\mu_f}\br).
\leqno(10)$$
We can prove this equality using a formula for Hodge ideals in the weighted homogeneous and non-degenerate cases in \cite{Zh} which can be extended easily to the semi-weighted-homogeneous case applying \cite[Theorem 0.9]{mos} (where $f$ is {\it not\1} assumed convenient), see (2.5) below. In the non-degenerate case, we need a rather non-trivial assertion \cite[Proposition B.1.2.3]{BGMM} for the proof of \cite[Theorem 5.5]{Zh}, where the assumption that $f$ is {\it convenient\1} is required. We can show, however, that this condition can be avoided, see Theorem~(2.7) below.
\sk
We have the following variants of Theorem~1.
\msn
{\bf Theorem~2.} {\it Assume $f\in\mm_{X,0}^3$, $\ep_f>0$, and $\mu_f\ne\tau_f$. Then $\Sp_f^{\HI}(t)\ne\Sp_f(t)$, and we have in the notation of $\1(3)$}
$$\al_{f,i}^{\HI}>\al_{f,\,\mu_f}\,\,\,(\forall\,i\notin I).
\leqno(11)$$
\msn
{\bf Theorem~3.} {\it Assume condition~{\rm (A)} is satisfied, and there is a monomial $g=\prod_ix_i^{\nu_i}\in\C\{x\}$ satisfying
$$fg\1\Om_f^n\ne 0\,\,\,(\h{that is,}\,\,\,fg\notin(\dd f))\q\h{and}\q\ga_f(g)+1>\al_{f,\,\mu_f}.
\leqno(12)$$
Then $\al_{f,\,\mu_f}^{\HI}>\al_{f,\,\mu_f}$, hence $\Sp_f^{\HI}(t)\ne\Sp_f(t)$.}
\ms
Their proofs are similar to Theorem~1 with $\ep_f>0$, see (3.2) below. 
In Theorem~2, we have $\ep_f>0$, that is, $\ga_f+1>\al_{f,\,\mu_f}$, if $\al_{f,n+1}+1>\al_{f,\,\mu_f}$. (Since the $\al_{f,i}$ are weakly increasing and $\dim\C\{x\}/\mm_{X,0}^2=n+1$, we have $\ga_f\ges\al_{f,n+1}$, where the equality holds if $w_{\max}\les 2\1w_{\min}$ in the semi-weighted-homogeneous case. The converse of the above assertion does not necessarily hold if $w_{\max}>2\1w_{\min}$.) Note that $\ga_f(g)$ is {\it not\1} defined unless condition~(A) is satisfied.
\sk
In the semi-weighted-homogeneous case with $f\in\mm_{X,0}^3$, it is known that the weights $w_i$ are unique and strictly smaller than $\tfrac{1}{2}$, see \cite{SaK}. If condition~(A) holds with $f\in\mm_{X,0}^3$, the following seems to be valid (as far as calculated):
$$\ga_f\les\tfrac{n+1}{3},\q\q\ep_f\les\tfrac{2n+1}{3}-(n{-}1)=\tfrac{4-n}{3}.
\leqno(13)$$
The last inequality would imply that $n\les 3$ when $\ep_f>0$.
\sk
In the {\it double point\1} case, that is, if $f\notin\mm_{X,0}^3$ (where $f\in\mm_{X,0}^2$ by the assumption that $Z$ has a singularity at 0), the weights of a weighted homogeneous polynomial are not necessarily unique (see \cite{SaK}), and $\ga_f-\al_{f,1}$ {\it can be arbitrary close to\1} 1 (consider for instance the case $f=x_1x_2+x_3x_4$). In this case, we have the following.
\msn
{\bf Proposition~1.} {\it Assume $f$ is a double point, that is, $f\notin\mm_{X,0}^3$, and $\mu_f\ne\tau_f$. Then $f\1\Om_f^n\subset V_{\HI}^{\,\al_{f,1}+2}$, hence $\al_{f,i}^{\HI}\ges\al_{f,1}+2$ for $i\notin I$. In particular, we have $\Sp_f^{\HI}(t)\ne\Sp_f(t)$ if $\al_{f,1}+2>\al_{f,\,\mu_f}$, that is, if $\al_{f,1}>\frac{n}{2}-1$.}
\msn
{\bf Proposition~2.} {\it Assume condition~{\rm (A)} holds, $f=h+x_n^2$ with $h\in\C\{x'\}:=\C\{x_1,\dots,x_{n-1}\}$, and moreover there is $g\in\C\{x'\}$ satisfying
$$fg\1\Om_f^n\ne 0\,\,\,(\h{that is,}\,\,\,fg\notin(\dd f))\q\h{and}\q v(g)+2>\al_{f,\,\mu_f},
\leqno(14)$$
where $v(g)$ is defined with $g$ viewed as an element of $\C\{x\}$.
Then $\al_{f,\,\mu_f}^{\HI}>\al_{f,\,\mu_f}$, hence $\Sp_f^{\HI}(t)\ne\Sp_f(t)$.}
\ms
In the double point case, the relation between Theorem~3 and Proposition~2 is not very clear. Indeed, if we have $f=u(h_2+h_1x_n+x_n^2)$ with $u$ invertible and $h_1,h_2\in\C\{x'\}$ after a coordinate change by the Weierstrass preparation theorem, we then get $f=u(h+x_n^2)$ with $h\in\C\{x'\}$ replacing $x_n$ with $x_n-h_1/2$. However, $f$ does not necessarily satisfy condition~(A) after these coordinate changes even if $f$ satisfies it before them.
\sk
Note finally that we may have $f^{n-1}\Om_f^n\ne 0$ in general. Extending the main theorems to this case is, however, rather difficult by the complexity of Hodge ideals, see Corollary~(2.6) and Example~(4.2) below. (The latter shows that the Hodge ideals are {\it not weakly decreasing even modulo the Jacobian ideal.})
\sk
In Section 1 we review some basics of the spectrum of isolated hypersurface singularities. In Section 2 we review Hodge ideals for $\Q$-divisors on complex manifolds. In Section 3 we prove the main theorems and propositions. In Section 4 we calculate some examples. In Appendix, we prove the key Proposition~(A.2) to the proof of Theorem~(2.7) after recalling some basics of Newton polyhedra.
\msn
{\bf Acknowledgement.} The first three authors are deeply grateful to Nero Budur for his valuable suggestion and encouragement. His suggestion was the starting point to study this problem. We thank the referee for useful comments to improve the paper.
\bs\bs
\vbox{\centerline{\bf 1. Spectrum}
\bsn
In this section we review some basics of the spectrum of isolated hypersurface singularities.}
\msn
{\bf 1.1.~Vanishing cohomology.} Let $f:(X,0)\to(\C,0)$ be a germ of a holomorphic function on a complex manifold having an isolated singularity at 0. Set $n=\dim X$. We have a canonical mixed Hodge structure on the vanishing cohomology $H^{n-1}(\Ff,\Q)$, where $\Ff$ denotes the Milnor fiber of $f$, see \cite{St2}. The $\la$-eigenspace of the semisimple part $T_s$ of the monodromy $T$ is denoted by $H^{n-1}(\Ff,\C)_{\la}$. (Note that $T$ is the {\it inverse\1} of the Milnor monodromy, see \cite{DS3}.)
\sk
The {\it spectrum} $\Sp_f(t)=\sum_{i=1}^{\mu_f}t^{\,\al_{f,i}}$ is a fractional power polynomial defined by
$$\aligned&\#\{i\mid\al_{f,i}=\be\}=\dim_{\C}\Gr_F^pH^{n-1}(\Ff,\C)_{\la}\\&\h{for}\q\be\in\Q,\,\,\,\,p=[n-\be],\,\,\,\la=\exp(-2\pi i\be),\endaligned
\leqno(1.1.1)$$
where $[\al]$ denotes the integer part of $\al\in\Q$, and the $\al_{f,i}$ are assumed weakly increasing, see \cite{St2}. These are positive rational numbers, and are called the {\it spectral numbers\1} or {\it exponents.} We have the symmetry
$$\al_{f,i}+\al_{f,j}=n\q\,\,\,\,(i+j=\mu_f+1).
\leqno(1.1.2)$$
This follows from the assertion that the weight filtration $W$ on $H^{n-1}(\Ff,\C)_{\la}$ coincides with the monodromy filtration which is associated with the action of $N:=\log T_u$, and is shifted by $n$ or $n{-}1$ depending on whether $\la=1$ or not (where $T_u$ is the unipotent part of the monodromy $T$), see {\it loc.\,cit.} This means that there are isomorphisms for $i>0:$
$$\aligned N^i:\Gr^W_{n-1+i}H^{n-1}(\Ff,\C)_{\la}&\simto\Gr^W_{n-1-i}H^{n-1}(\Ff,\C)_{\la}\q(\la\ne 1),\\ N^i:\Gr^W_{n+i}H^{n-1}(\Ff,\C)_1&\simto\Gr^W_{n-i}H^{n-1}(\Ff,\C)_1.\endaligned
\leqno(1.1.3)$$
Note that the $N$-{\it primitive part\1} can be defined by the kernel of $N^{i+1}$ ($i\ges 0$) so that the $N$-{\it primitive decomposition\1} holds.
\sk
It is known (see for instance \cite[4.11]{DS2}) that the {\it multiplicity\1} of the minimal (or maximal) exponent is 1, that is,
$$\#\bl\{i\in[1,\mu_f]\,\,\big|\,\,\al_{f,i}=\al_{f,1}\br\}=\#\bl\{i\in[1,\mu_f]\,\,\big|\,\,\al_{f,i}=\al_{f,\,\mu_f}\br\}=1.
\leqno(1.1.4)$$
\msn
{\bf 1.2.~Brieskorn lattices.} In the notation of (1.1), we have the {\it Brieskorn lattice\1} (see \cite{Br})
$$H''_f:=\Om_{X,0}^n/\ddd f{\wedge}\ddd\Om_{X,0}^{n-2}.$$
This is a free module of rank $\mu_f$ over $\C\{t\}$ and also over $\C\{\!\{\dd_t^{-1}\}\!\}$, where the action of $\dd_t^{-1}$ is defined by $\dd_t^{-1}[\om]=[\ddd f\wedge\eta]$ for $\ddd\eta=\om$, see \cite{SS}, \cite{bl}. Recall that $\C\{\!\{\dd_t^{-1}\}\!\}$ is defined by
$$\bl\{\1\msum_{k\in\N}\,c_k\1 \dd_t^{-k}\,\big|\,\msum_{k\in\N}\,c_k\1 x^k/k!\in\C\{x\}\1\br\}.$$
\sk
The localization $G_f:=H''_f[\dd_t]$ by the action of $\dd_t^{-1}$ is called the {\it Gauss-Manin system}. It is a free $\C\{\!\{\dd_t^{-1}\}\!\}[\dd_t]$-module of rank $\mu_f$ with an action of $t$, and is a regular holonomic $\D_{X,0}$-module. So it has the $V$-filtration of Kashiwara \cite{Ka2} and Malgrange \cite{Ma2} indexed by $\Q$ and such that $\dd_tt-\be$ is nilpotent on $\Gr_V^{\be}G_f$. The filtration $V$ on $H''_f$ is the induced filtration of the $V$-filtration by the inclusion $H''_f\into G_f$. (Note that the filtration $V$ in \h{\it loc.\,cit.} was indexed by $\Z$, instead of $\Q$. The $V$-filtration indexed by $\Q$ was influenced by the theory of asymptotic Hodge filtration \cite{Va1} using the asymptotic expansions of period integrals.)
\sk
On the other hand, there are canonical surjection
$$H''_f\,\onto\,\Om_f^n:=\Om_{X,0}^n/\ddd f{\wedge}\Om_{X,0}^{n-1},
\leqno(1.2.1)$$
and the filtration $V$ on $\Om_f^n$ is defined as the quotient filtration of the $V$-filtration on $H''_f$.
\sk
By \cite{SS} (see also \cite{Va1}) we have the canonical isomorphisms
$$\Gr_V^{\be}\Om_f^n=\Gr_F^pH^{n-1}(\Ff,\C)_{\la}\q(p=[n-\be],\,\,\la=e^{-2\pi i\be}).
\leqno(1.2.2)$$
This assertion is shown by using the Hodge filtration on $G_f$ defined by
$$F^pG_f:=\dd_t^{n-1-p}H''_f\q(p\in\Z).
\leqno(1.2.3)$$
since it induces the Hodge filtration on $H^{n-1}(\Ff,\C)$ via the isomorphisms
$$\Gr_V^{\al}G_f=H^{n-1}(\Ff,\C)_{\la}\q(\al\in(0,1],\,\la=e^{-2\pi i\al}).
\leqno(1.2.4)$$
\msn
{\bf Remark~1.2}\,(i). If we replace the Hodge filtration $F$ with the asymptotic Hodge filtration in \cite{Va1} which can be defined by using
$$t^{p-n+1}H''_f\subset H''_f[t^{-1}],
\leqno(1.2.5)$$
instead of (1.2.3), then the isomorphisms in (1.2.2) hold after taking $\Gr^W_i$.
\msn
{\bf Remark~1.2}\,(ii). We have the inclusions
$$f\,V^{\be}\Om_f^n\subset V^{\be+1}\Om_f^n\q\q(\be\in\Q),
\leqno(1.2.6)$$
since $f$ is identified with $\Gr_Ft$ on $\Gr_F^0G_f=\Om_f^n$.
\msn
{\bf 1.3.~Microlocal $V$-filtration.} Let $(i_f)_+\OO_X$ be the direct image of $\OO_X$ as a $\D$-module by the graph embedding $i_f:X\to X\times\C$ associated with $f$. This is identified with
$$M_f:=\OO_X[\dd_t]\delta(t{-}f),$$
where $t$ is the coordinate of $\C$. The action of $\dd_{x_i}$ on $M_f$ is given by
$$\dd_{x_i}\bl(g\,\dd_t^k\delta(t{-}f)\br)=(\dd_{x_i}g)\dd_t^k\delta(t{-}f)-g(\dd_{x_i}f)\dd_t^{k+1}\delta(t{-}f),
\leqno(1.3.1)$$
where $g\in\OO_X$. We have the $V$-filtration of Kashiwara \cite{Ka2} and Malgrange \cite{Ma2} on $M_f$ indexed by $\Q$ such that $\dd_tt-\be$ is nilpotent on $\Gr_V^{\be}M_f$.
\sk
We have the {\it algebraic microlocalization}
$$\M_f:=\OO_X[\dd_t,\dd_t^{-1}]\delta(t{-}f).
\leqno(1.3.2)$$
It has the Hodge filtration $F$ by the order of $\dd_t^{-1}$ together with the $V$-filtration defined by using the $V$-filtration on $M_f$, see \cite{mic}. The {\it $V$-filtration} and the {\it microlocal $V$-filtration\1} on $\OO_X$ are defined respectively by the filtered isomorphisms (see \cite{mic}):
$$(\OO_X,V)=\Gr_F^0(M_f,V),\q(\OO_X,\V)=\Gr_F^0(\M_f,V).
\leqno(1.3.3)$$
\msn
{\bf Remark~1.3}\,(i). By the construction (see \h{\it loc.\,cit.} and also (2.5.4) below), we have
$$V^{\al}\OO_X=\V^{\al}\OO_X\q(\al\les 1).
\leqno(1.3.4)$$
and these can be identified essentially with the {\it multiplier ideals} (see \cite{BS1}) except for the difference in the index of filtration, where the $V$-filtrations are left-continuous, but the multiplier ideals are right continuous, see also \cite{MSS2}.
\msn
{\bf Remark~1.3}\,(ii). The {\it Bernstein-Sato polynomial} $b_f(s)$ is defined as the minimal polynomial of the action of $s:=-\dd_tt$ on
$$N_f:=\D_X[s]f^s\big/\D_X[s]f^{s+1},$$
where $f^s$ is identified with $\delta(t{-}f)$ (see (2.2) below) so that we have the inclusion
$$\D_X[s]f^s=\D_X\langle s,t\rangle f^s\subset M_f.$$
The {\it microlocal Bernstein-Sato polynomial} $\bt_f(s)$ is defined as the minimal polynomial of the action of $s:=-\dd_tt$ on
$$\widetilde{N}_f:=\D_X\langle s,\dd_t^{-1}\rangle\,\delta(t{-}f)\big/\D_X\langle s,\dd_t^{-1}\rangle\,\dd_t^{-1}\delta(t{-}f),$$
and coincides with the {\it reduced Bernstein-Sato polynomial} $b_f(s)/(s{+}1)$, see \cite{mic}.
\sk
We have the inclusions
$$\Gr_V^{\al}\OO_X\into\Gr_V^{\al}N_f\q(\al<\al_f+1),
\leqno(1.3.5)$$
\vskip-6mm
$$\Gr_{\V}^{\al}\OO_X\into\Gr_V^{\al}\widetilde{N}_f\q(\al<\widetilde{\al}_f+1),
\leqno(1.3.6)$$
where $\al_f$, $\widetilde{\al}_f$ are, up to sign, the maximal root of $b_f(s)$, $\bt_f(s)$ respectively, and we have $\al_f=\min(\widetilde{\al}_f,1)$. For the proof of (1.3.6), we need the inclusion
$$F^1\M_f\subset V^{\alt_f+1}\M_f,$$
since the definition of $\V$ uses the graded quotient as in the second isomorphism of (1.3.3), and {\it not the inclusion\1} $\OO_X\delta(t{-}f)\into\M_f$. 
\sk
As for (1.3.5), we cannot get any information for $\al\in(1,\al_f+1)$, since the source vanishes for such $\al$ by the {\it periodicity\1} of the $V$-filtration (or multiplier ideals)
$$f\1V^{\al}\OO_X=V^{\al+1}\OO_X\q(\al>0).
\leqno(1.3.7)$$
\sk
From (1.3.5--6) one can deduce some relation between the (microlocal) jumping coefficients and the roots of (microlocal) Bernstein-Sato polynomial, since the action of $s+\al$ is nilpotent on $\Gr_V^{\al}N_f$, $\Gr_V^{\al}\widetilde{N}_f$. For instance, (1.3.6) implies that any microlocal jumping coefficient (see \cite{MSS2}) strictly smaller than $\alt_f+1$ is a root of the microlocal Bernstein-Sato polynomial up to sign. A similar assertion for usual jumping coefficients is well-known, see \cite{ELSV}.
\msn
{\bf Remark~1.3}\,(iii). It is possible to show the last assertion in Remark~(1.3)\,(ii) above by using $|f^2|^s$ and integration by parts, or equivalently, derivation as a distribution, although certain delta functions or locally non-integrable functions (that is, not belonging to $L^1_{\rm loc}(X)$) might appear in the intermediate stages, when one applies the derivation as a distribution repeatedly. Here it is highly recommended to use analytic continuation in $s$ so that some equality can be shown for ${\rm Re}\,s\gg 0$ forgetting the intermediate stages. Sometimes this point does not seem to be explained sufficiently to the reader.
\msn
{\bf Remark~1.3}\,(iv). The minimal exponents $\al_f,\alt_f$ in (1.3.5--6) are invariant under a {\it non-characteristic restriction} to a closed submanifold $X'\subset X$, for instance, if $X'$ transversally intersects any strata of a Whitney stratification compatible with $\psi_f\C_X$.
This follows for instance from \cite[Theorem~1.1 and Lemma~4.1]{DMST}, since the minimal exponents can be determined by the Hodge filtration on the $\Gr_V^{\al}M_f$ using \cite[2.1.4 and 2.2.3]{mic}.
\sk
One can use this invariance to deduce the semi-continuity assertion in Theorem~E\,(2) from Theorem~E\,(1) in \cite{MP3} by induction on the dimension of the parameter space $T$. Indeed, there is a non-empty Zariski open subset $T'\subset T$ such that the fiber of $t\in T'$ is transversal to any strata of the above stratification, and the function $x\mapsto\alt_{D,x}$ is semi-continuous.
\msn
{\bf 1.4.~Coincidence of the two quotient $V$-filtrations.} In this subsection we prove the following (see also \cite[(4.11.4)]{DS2}).
\msn
{\bf Proposition~1.4.} {\it The quotient $V$-filtration on $\Om_f^n\cong\C\{x\}/(\dd f)$ in $(1.2)$ defined by using the surjection $(1.2.1)$ coincides with the quotient $V$-filtration induced from the microlocal $V$-filtration on $\OO_{X,0}=\C\{x\}$ using the canonical surjection $\C\{x\}\,\onto\,\C\{x\}/(\dd f)$.}
\msn
{\it Proof.} This is proved by using the {\it strictness\1} of the bifiltered de Rham complex
$$\bl({\rm DR}(\M_{f,0});F,V\br).$$
Note first that the highest cohomology
$$H^n\bl({\rm DR}(\M_{f,0})\br)$$
coincides with the {\it Gauss-Manin system\1}
$$H^n\bl({\rm DR}(M_{f,0})\br)=G_f,$$
since the action of $\dd_t$ on the latter is bijective. The filtrations on the top cohomology induced by $V$, $F$ coincide respectively with the $V$-filtration of the Gauss-Manin system $G_f$ (see \cite[(4.11.5)]{DS2}) and the filtration defined by the shifted Brieskorn lattices $\dd_t^{-p}H''_f$ ($p\in\Z$).
\sk
The strictness of the above bifiltered complex means the {\it bistrictness\1} for the two filtrations $F$, $V$, that is, the three filtrations $F,V,G$ on each component form compatible filtrations, where $G$ is defined by the kernel and image of the differential of the complex. Then $\Gr_F^{\ssb}$ and $H^{\ssb}$ commute in a compatible way with the filtration $V$, see also \cite[Corollary~1.2.13]{mhp}. Here $\Gr_F^{\ssb}$ of the de Rham complex is identified with the {\it Koszul complex}
$$K^{\ssb}(\C\{x\};f_1,\dots,f_n)$$
associated with the multiplications by the partial derivatives $f_i:=\dd_{x_i}f$ ($i\in[1,n]$) on $\C\{x\}$. The filtration induced by $V$ on
$$\Gr_F^0(\M_{f,0})=\OO_{X,0}=\C\{x\}$$
coincides with the microlocal $V$-filtration on $\C\{x\}$ by definition, see (1.3.3).
\sk
These imply that the quotient filtration $V$ on $\C\{x\}/(\dd f)$ induced by the microlocal $V$-filtration on $\C\{x\}$ coincides with the one induced by the surjection (1.2.1), since the canonical filtered isomorphism
$$H^n\Gr_F^0({\rm DR}(\M_{f,0}),V)=\Gr_F^0H^n({\rm DR}(\M_{f,0}),V)
\leqno(1.4.1)$$
provides the filtered isomorphism
$$H^n(K^{\ssb}(\C\{x\};f_1,\dots,f_n),V)=\Gr_F^0(G_f,V).
\leqno(1.4.2)$$
So Proposition~(1.4) follows.
\msn
{\bf Remark~1.4.} In the weighted homogeneous isolated singularity case, Proposition~(1.4) follows also from \cite[Proposition in 2.2]{hi} where an explicit description of the microlocal $V$-filtration is given for weighted homogeneous isolated singularities. (Here the assumption that any weight $w_i$ is the inverse of an integer is not necessary by using \cite[Remark (ii) in Section 4.1]{mos} as is remarked by Mingyi Zhang.)
\msn
{\bf 1.5.~Weighted homogeneous case.} There is a well-known formula for the spectrum of a weighted homogeneous polynomial with an isolated singularity
$$\Sp_f(t)=\mprod_{i=1}^n\,(t^{w_i}-t)/(1-t^{w_i}),
\leqno(1.5.1)$$
where the $w_i$ are the weights of variables $x_i$ associated with $f$. This formula is conjectured in \cite{St2}, and follows essentially from \cite{St1} where a generalization of the Griffiths theorem to the Milnor cohomology of weighted homogeneous isolated hypersurface singularities is proved (although the relation with the monodromy should be added).
\sk
Actually the equality (1.5.1) follows also from \cite{SS}, \cite{Va1}. Indeed, by an argument similar to the proof of \cite[Proposition~3.3]{exp}, we can show that
$$\dd_tt[x^{\nu}\ddd x]=\ell_w(\nu)[x^{\nu}\ddd x]\q\h{with}\q\ell_w(\nu):=\msum_i\,w_i(\nu_i+1),
\leqno(1.5.2)$$
where $x^{\nu}=\prod_ix_i^{\nu_i}$, and $\ddd x:=\ddd x_1\wedge\cdots\wedge\ddd x_n$, see also \cite[Section 2.2]{DS3}. Define a filtration $V_w$ on $\C\{x\}$ by
$$V_w^{\be}\C\{x\}:=\msum_{\ell_w(\nu)\ges\be}\,\C\{x\}x^{\nu}.
\leqno(1.5.3)$$
Then (1.5.2) implies the following:
$$\h{The $V$-filtration on $H''_f$ is given by the quotient filtration of $V_w$.}
\leqno(1.5.4)$$
This and (1.2.2) imply (1.5.1), since we have the equality
$${\rm HP}\bl(\C\{x\}/(\dd f),V_w\br)=\mprod_{i=1}^n\,(t^{w_i}-t)/(1-t^{w_i}),
\leqno(1.5.5)$$
where the left-hand side denotes the Hilbert Poicar\'e series of the filtered vector space $\bl(\C\{x\}/(\dd f),V_w\br)$, that is, the one for the graded vector space $\Gr_{V_w}^{\ssb}\bl(\C\{x\}/(\dd f)\br)$.
\msn
{\bf Remark~1.5}\,(i). In \cite{JKY} the equality (1.5.5) is shown by using an inductive argument. However, there are shorter methods as follows:
\sk
(a) As in \cite[(4.2.2)]{bl}, \cite[(4.1.1)]{mos}, consider the flat morphism
$$\rho:=(f_1,\dots,f_n):\C^n\to\C^n,$$
so that $\C[x]$ is a {\it free graded module\1} over $\C[y]$ generated freely by a $\C$-basis of $\C[x]/(\dd f)$ (which is isomorphic to $\C\{x\}/(\dd f))$ using the graded version of Nakayama's lemma, where $y_1,\dots,y_n$ are the coordinates of the target space $\C^n$ so that $\rho^*y_i=f_i\,(:=\dd f/\dd x_i)$. This implies the equality
$$\mprod_{i=1}^n\,(1-t^{w_i})^{-1}=t^{-\al_1}\1\Sp_f(t)\1\mprod_{i=1}^n\,(1-t^{1-w_i})^{-1}\q\h{in}\,\,\,\Z[[t]],
\leqno(1.5.6)$$
since $\deg f_i=1-w_i$ and $t^{-\al_1}\1\Sp_f(t)$ is the Hilbert-Poincar\'e series of $\C[x]/(\dd f)$ with $V$-filtration {\it not shifted by\1} $\al_{f,1}=\sum_iw_i$ as in (1.5.2--3). Indeed, the left-hand-side of (1.5.6) is the Hilbert-Poincar\'e series of the graded vector space $\C[x]$ with $\deg x_i=w_i$, and similarly for $\C[y]$ with $\deg y_i=1-w_i$.
\ms
(b) Observe first that the coefficient of $X^p$ in the polynomial
$$\mprod_{i=1}^n(X+t^{w_i})(1-t^{w_i})^{-1}\in\Q[[t]][X]
\leqno(1.5.7)$$
is the Hilbert-Poincar\'e series of the graded vector space of algebraic $(n{-}p)$-forms $\Om^{n-p}$, where $\deg\ddd x_i=w_i$. Since the Koszul complex $(\Om^{\ssb},\ddd f\wedge)$ gives a graded free resolution of $\Om_f^n$, and $f$ and $\ddd f$ have degree 1 by assumption, it is then enough to substitute $-t$ in $X$ to get the Hilbert-Poincar\'e series of $\Om_f^n$. Here the minus sign is needed to calculate the Euler characteristic of the Koszul complex (up to a sign).
\msn
{\bf Remark~1.5}\,(ii). If $f$ is a weighted homogenous polynomial with an isolated singularity, the roots of the reduced Bernstein-Sato polynomial have multiplicity 1, and coincide with the spectral numbers up to sign. Admitting (1.5.1), (1.5.5), this was proved in \cite{Ma1} or \cite{Sat}. (In the latter, an argument similar to \cite[Remark 4.2(i)]{mos} was used.)
\msn
{\bf 1.6.~Semi-weighted-homogeneous case.} A holomorphic function $f\in\C\{x\}$ having an {\it isolated singularity\1} at 0 is called {\it semi-weighted-homogeneous with weights} $w_1,\dots,w_n>0$ if there is a decomposition $f=f_1+f_{>1}$ such that $f_1$ a weighted homogeneous polynomial of weights $w_i$ having an isolated singularity at the origin, and $f_1$, $f_{>1}$ are $\C$-linear combinations of monomials $x^{\nu}:=\prod_ix_i^{\nu_i}$ satisfying $\msum_i\,w_i\nu_i=1$ and $\msum_i\,w_i\nu_i>1$ respectively, where $\nu=(\nu_1,\dots,\nu_n)\in\N^n$. This is equivalent to that $f$ is a $\mu$-constant deformation of a weighted homogeneous polynomial with an isolated singularity, see \cite{Va3}.
\sk
Since the spectrum is invariant by $\mu$-constant deformations (see for instance \cite{Va2}), we have
$$\Sp_f(t)=\Sp_{f_1}(t)=\mprod_{i=1}^n\,(t^{w_i}-t)/(1-t^{w_i}).
\leqno(1.6.1)$$
\sk
For $g=\msum_{\nu}c_{\nu}x^{\nu}\in\C\{x\}$, define
$$v(g):=\min\bl\{\msum_i\,w_i(\nu_i+1)\mid c_{\nu}\ne 0\br\}.
\leqno(1.6.2)$$
It is known that (1.5.4) also holds in the semi-weighted-homogeneous case by defining $V_w$ as in (1.5.3) (using a 1-parameter $\mu$-constant deformation where we get a constant family of Gauss-Manin systems, and taking the completion by the filtration $V_w$). We have
$$v(1)=\al_{f,1}=\msum_{i=1}^n\,w_i.
\leqno(1.6.3)$$
Put
$$\ga_f(g):=\max\bl\{v(x_i\1g)\}_{i\in[1,n]},\q\q\ga_f:=\ga_f(1).
\leqno(1.6.4)$$
We have
$$\ga_f=\al_{f,1}+w_{\max}\q\h{with}\q w_{\max}:=\max\{w_i\}_{i\in[1,n]}.
\leqno(1.6.5)$$
\msn
{\bf 1.7.~Non-degenerate Newton boundary case.} Let $\Gp(f)$ be the {\it Newton polyhedron\1} of $f\in\C\{x\}$ at $0$, that is, the convex hull of the union of $\nu+\R_{\ges 0}^n$ for $\nu\in\Spx\1f$ with
$$\Spx\1f:=\{\nu\in\N^n\mid a_{\nu}\ne 0\}\q\h{if}\q f=\msum_{\nu}\,a_{\nu}x^{\nu}\in\C\{x\}.
\leqno(1.7.1)$$
We will denote by $\Ga(f)$ the convex hull of $\Spx f$.
We say that $f$ has {\it non-degenerate Newton boundary} at the origin, if for any {\it compact\1} face $\si\subset\Gp(f)$, we have
$$\bl\{x_i\dd_{x_i}f_{\si}=0\,\,\,(\forall\,i\in[1,n])\br\}\cap(\C^*)^n=\emptyset,
\leqno(1.7.2)$$
where $f_{\si}:=\mopl_{\nu\in\si}\,a_{\nu}x^{\nu}$, see \cite{Ko}. (This definition is equivalent to the one in \cite{exp}, see \cite{Ko}.)
\sk
For $g\in\C\{x\}$, set
$$v(g):=\max\bl\{a\in\R\mid\ob^n\,{+}\,\Spx\1g\,\subset\,a\,\Gp(f)\br\},
\leqno(1.7.3)$$
with $\ob^n:=(1,\dots,1)$. Define $\ga_f(g)$, $\ga_f$ by (1.6.4). Note that $\al_{f,1}=v(1)$, and the $V$-filtration on $\Om_f^n\cong\C\{x\}/(\dd f)$ in (1) is induced by the Newton filtration $V_N^{\ssb}$ on $\C\{x\}$ defined by
$$V_N^{\be}\C\{x\}:=\{v(g)\ges\be\}\,\,\,\h{(see \cite{exp}, \cite{VK})}.
\leqno(1.7.4)$$
(Here the convenience condition is not necessary in the isolated singularity case as is shown in Corollary~A in Appendix below.)
\sk
We say that $f$ is {\it convenient\1}, if the intersection of $\Gp(f)$ with every coordinate axis is non-empty. In the isolated hypersurface singularity case, we may assume $f$ is convenient by adding monomials $x_i^{a_i}$ to $f$ for $a_i\gg 0$ if necessary. This is allowed by the {\it finite determinacy\1} of holomorphic functions with isolated singularities as is used in \cite{Br}, \cite{SS}, \cite{Va1}, etc. It is expected that $\ga_f(g)$ remains unchanged by adding $x_i^{a_i}$ to $f$ if $a_i\gg 0$ (depending on $g$).
\msn
{\bf Remark~1.7.} Condition~(1.7.2) is equivalent to that the subvariety of $(\C^*)^n$ defined by $f_{\si}$ is reduced and smooth. Indeed, any {\it compact\1} face $\si\subset\Gp(f)$ is contained in a hyperplane $H_{\ell}\subset\R^n$ defined by a linear equation $\ell(\nu):=\sum_{i=1}^nc_i\nu_i-c_0=0$ with $c_i$ {\it strictly positive} $(\forall\,i\in[0,n])$. Note that the non-vanishing of $c_0$ is required to claim that $f_{\si}$ is a linear combination of the $x_i\dd_{x_i}f_{\si}$. (This equivalence does not hold, that is, the smoothness of $f_{\si}^{-1}(0)$ is not enough, in the case one considers a Newton polyhedron at infinity associated with a {\it Laurent polynomial\1} $f$ unless $\si$ is contained in a hyperplane $H_{\ell}$ with $c_0\ne 0$.)
\msn
{\bf 1.8.~Spectrum for the general hypersurface case.} Let $f:(X,0)\to(\C,0)$ be a non-constant holomorphic function (which does not necessarily have an isolated singularity). We define the {\it Steenbrink spectrum\1} by
$$\Sp_f(t):=\msum_{j=0}^{n-1}\,(-1)^j\Sp_f^j(t),
\leqno(1.8.1)$$
where $\Sp^j_f(t)$ is the $j$\,th spectrum defined by replacing $H^{n-1}(\Ff,\C)$ with $\Ht^{n-1-j}(\Ff,\C)$ in (1.1.1) for $j\in[0,n-1]$ ($\Ht$ denotes the reduced cohomology). This is a fractional power polynomial as in the isolated singularity case (although the definition in \cite{St3} is slightly different by the division by $t$).
\sk
It is sometimes better to use
$$\Sp'_f(t):=t^n\,\Sp_f\bl(\tfrac{1}{t}\br),
\leqno(1.8.2)$$
for instance, in the proof of the Steenbrink conjecture, see \cite[Section 2.2]{ste}. Note that $\Sp'_f(t)$ is called the {\it Hodge spectrum\1} by Denef and Loeser in the definition before \cite[Corollary 6.24]{DL}, although $\Sp_f(t)$ is so called in \cite{Bu}. These can be different, since the symmetry (1.1.2) as in the isolated singularity case does not necessarily hold. (For instance, the jumping coefficients in $(0,1)$ smaller than the minimal exponents at $x\in X\setminus\{0\}$ are Steenbrink spectral numbers in certain cases, but this is not clear for the Hodge spectrum since there may be a cancellation among the $\Sp^{\prime j}_f(t)$.)
\sk
The difference between $\Sp_f(t)$ and $\Sp'_f(t)$ is closely related to the one between $i_0^*$ and $i_0^!$ (where $i_0:\{0\}\into X$ is the inclusion), since the vanishing cycle sheaves are {\it self-dual\1} up to a Tate twist depending on the monodromy eigenvalue, see \cite[(2.6.2)]{mhm}. Note that the mixed Hodge structure on the vanishing cohomology in the non-isolated singularity case is defined by applying the cohomological functor $H^ki_0^*$ to the vanishing cycle Hodge module (which is denoted by $\varphi_f\Q_{h,X}[n{-}1]$ in this paper).
\bs\bs
\vbox{\centerline{\bf 2. Hodge ideals on complex manifolds}
\bsn
In this section we review Hodge ideals for $\Q$-divisors on complex manifolds.}
\bsn
{\bf 2.1.~Hodge ideals in the analytic case.} The Hodge ideals are defined in the algebraic case, see \cite{MP1}, \cite{MP2}, \cite{MP3}. However, they can be defined also in the analytic setting, since so are Hodge modules (although the standard functors between the derived categories cannot be defined in general). This may be useful even in the isolated singularity case. Indeed, although analytic isolated hypersurface singularities are always algebraizable, the independence of algebraization for Hodge ideals does not seem to be completely trivial unless these are analytically defined.
\ms
Let $D=\msum_k\,\al_kZ_k$ ($\al_k\in\Q_{>0}$) be an effective $\Q$-divisor on a complex manifold (or a smooth complex algebraic variety) $X$ with $Z_k$ the reduced irreducible components. Set
$$Z:=\mcup_k\,Z_k,\q U:=X\setminus Z\q\h{with $\,j_U:U\into X\,$ the inclusion.}$$
Locally on $X$, there is a unique $\C$-local system $L_D$ of rank 1 on $U$ whose local monodromy around $Z_k\setminus{\rm Sing}\,Z$ is given by the multiplication by $e^{2\pi i\al_k}$ for any $k$. This is unique by considering the tensor product of one local system with the dual of another. Indeed, for a contractible open subset $B\subset X$, we have the vanishing of the first homology group $H_1(B\setminus{\rm Sing}\,Z,\Z)$ which is the abelianization of the fundamental group (since ${\rm Sing}\,Z\subset X$ has codimension at least 2).
Its existence can be shown by considering a $\C$-local system of rank 1 on $\C^*$ whose monodromy is given by the multiplication by $e^{2\pi i\al}$ ($\al\in\Q_{>0}$), and taking its pull-back by a holomorphic function $f$ such that
$$D=\al\1({\rm div}\1 f)\q\h{(hence}\,\,\,\,\al^{-1}\al_k\in\Z_{>0}).
\leqno(2.1.1)$$
Here we may assume GCD$(\al^{-1}\al_k)=1$ (although this condition is unstable by shrinking $X$).
\sk
Let $\Lc_D(*Z)$ be the meromorphic extension of $\OO_U\otimes_{\C}L_D$. This is a regular holonomic $\D_X$-module corresponding to $\RR(j_U)_*L_D$, and can be defined also by the direct image of the meromorphic extension on an embedded resolution of $Z\subset X$ (see \cite{De}). It has a canonical meromorphic connection together with a locally free $\OO_X$-submodule $\Lc_D$ such that
$$\Lc_D|_U=\OO_U\otimes_{\C}L_D,
\leqno(2.1.2)$$
and the restriction of $\Lc_D$ to $U':=X\setminus{\rm Sing}\,Z$ is stable by the logarithmic connection with eigenvalue of its residue along $Z_k\setminus{\rm Sing}\,Z$ equal to $-\al_k$.
This is unique, since so is its restriction to $U'$ ({\it loc.\,cit.}), and the Hartogs theorem implies that
$$\Lc_D=(j_{U'})_*j_{U'}^{-1}\Lc_D,
\leqno(2.1.3)$$
with $j_{U'}:U'\into X$ the inclusion. Its existence can be shown by using a holomorphic function $f$ as above, since we have locally on $X$ the following isomorphism of regular holonomic $\D_X$-modules (see also \cite{rp}):
$$\Lc_D(*Z)=\OO_X(*Z)f^{-\al}.
\leqno(2.1.4)$$
\sk
We can show that $\Lc_D(*Z)$ is isomorphic to a direct factor of the underlying $\D_X$-module of a mixed Hodge module by using the action of the covering transformation group of a ramified covering defined by $\widetilde{X}:=\{z^e=f\}\subset X\times\C$ with $e$ a positive integer satisfying $e\al\in\N$.
Here we first take the open direct image of the constant Hodge module on $\widetilde{X}\setminus\{z=0\}$ under the open inclusion into $\widetilde{X}$ (see \cite[Proposition 2.11]{mhm}), and then the cohomological direct image under the projection $X\times\PP^1\to X$ (see for instance \cite[Proposition 2.14]{mhm}). For the algebraic case, see also \cite{MP2}. (It seems easier to show that the Hodge ideals are independent of a choice of $f$ in the analytic case; for instance, we can take $u^{1/e}$ for an invertible function $u\in\OO_{X,0}$, although we would need an \'etale neighborhood in the algebraic case.)
\sk
We then get the Hodge filtration $F$ on $\Lc_D(*Z)$ with the inclusions of $\OO_X$-modules
$$F_p\bl(\Lc_D(*Z)\br)\subset\Lc_D(p\1Z)\subset\Lc_D(*Z)\q(p\in\N),
\leqno(2.1.5)$$
where $\Lc_D(pZ)$ denotes the canonical image of $\Lc_D\otimes_{\OO_X}\OO_X(pZ)$ in $\Lc_D(*Z)$.
Note that this Hodge filtration depends only on the $\D$-module $\Lc_D(*Z)$ by the formula for the Hodge filtration on the open direct images (see \cite[3.2.3.2]{mhp}) applied to $j_U:U\into X$.
\sk
The Hodge ideals $I_p(\al Z)\subset\OO_X$ are then defined by the equality of $\OO_X$-submodules
$$F_p\bl(\Lc_D(*Z)\br)=I_p(\al Z)\,\Lc_D(pZ)\q\h{in}\,\,\,\Lc_D(pZ)\q(p\in\N).
\leqno(2.1.6)$$
Here the right-hand side is the $\OO_X$-submodule of $\Lc_D(pZ)$ consisting of $hv$ with $h\in I_p(\al Z)$ and $v$ {\it any\1} local generator of $\Lc_D(pZ)$.
\sk
By definition the Hodge ideals {\it modulo multiplications by locally principal ideals of} $\OO_X$ depend only on the $\al_k$ mod $\Z$ (since so is $\Lc_D$), see also \cite[Lemma 4.4]{MP2}. So we may assume as in \cite[Theorem A]{MP3}
$$\al_k\in(0,1]\,\,\,(\forall\,k),\q\h{that is,}\q\lceil D\rceil=Z.
\leqno(2.1.7)$$
The last condition is equivalent to that $\Lc_D|_{U'}$ is the Deligne extension with eigenvalues of the residues contained in $[-1,0)$, where $U'$ is as in (2.1.3). Set as in \cite{MP3}
$$I''_p(D):=I_p(D)\,\1\OO_X\bl(p(Z{-}{\rm div}\1f)\br)\subset\OO_X,
\leqno(2.1.8)$$
where the right-hand side is the product of the two ideals of $\OO_X$. This means that
$$F_p\bl(\Lc_D(*Z)\br)=I''_p(\al Z)\,\Lc_D\bl(p({\rm div}\1f)\br)\q\h{in}\,\,\,\Lc_D(*Z)\q(p\in\N).
\leqno(2.1.9)$$
If the $\al_k$ are {\it independent of\1} $k$, we assume $\al_k=\al$ ($\forall\,k$) so that $D=\al Z$ with $Z={\rm div}\1f$, and hence $I''_p(D)=I_p(D)$.
\msn
{\bf Remark~2.1}\,(i). In the {\it irreducible\1} divisor case with $D=\al Z$ ($\al>0$), one could define the ({\it microlocal\1$)$ Hodge ideals\1} by
$$I(\al Z):=I_p(\al' Z)\q\h{with}\q\al'\in(0,1],\,\,\,p=\al-\al'\in\N.
\leqno(2.1.10)$$
This is, however, impossible in the {\it reducible\1} case unless the $\lceil\al_k\rceil$ are {\it independent of $k$.}
\msn
{\bf Remark~2.1}\,(ii). If $Z$ is smooth (in particular, reduced), then the Hodge ideals $I(\al Z)$ coincide with $\OO_X$. Indeed, the corresponding assertion for the microlocal $V$-filtration is well-known (since the vanishing cycles vanish), and we can use the coincidence modulo $(f)$ (hence modulo the maximal ideal of $\OO_{X,x}$ for $x\in Z$). This implies that the support of the quotient $\OO_X/I(\al Z)$ is contained in the singular locus of $Z$.
\ms
Assuming (2.1.7) and choosing a holomorphic function $f$ as in (2.1.1), we can describe the Hodge ideals more explicitly using the {\it two isomorphisms\1} explained in the following two subsections.
\msn
{\bf 2.2.~The first isomorphism.} We denote by $i_f:X\into Y:=X\times\C$ the graph embedding with $t$ the coordinate of $\C$. Let $M$ be a holonomic $\D_X$-module such that the action of $f$ is {\it bijective}. It is well-known to specialists that there is a natural isomorphism of $\D_Y$-modules
$$(i_f)_+M\,\bl(=M[\dd_t]\1\delta(t{-}f)\br)\,\simto\,M[s]f^s.
\leqno(2.2.1)$$
This is noted essentially in \cite[p.\,110]{Ma1} in the case $M=\OO_X(*Z)$. The right-hand side of (2.2.1) has a structure of a left $\D_X\langle s,t\rangle$-module such that the action of $t$ is bijective. (We use $\langle s,t\rangle$ rather than $[s,t]$, since $s$, $t$ do not commute.) Here the action of $t$ on $M[s]f^s$ is defined by
$$t(m\1s^jf^s)=(fm)(s{+}1)^jf^s\q\q(m\in M),
\leqno(2.2.2)$$
(as is written in \cite{Ka1} in the case $M=\OO_X$). The action of $\OO_X$ on $M[s]f^s$ is a natural one, but the action of a vector field $\xi$ on $X$ is twisted as follows:
$$\xi(m\,s^jf^s)=\xi(m)\1s^jf^s+(\xi(f)/f)m\,s^{j+1}f^s\q\q(m\in M).
\leqno(2.2.3)$$
The $\D_Y$-module structure on $M[s]f^s$ together with the $\D_Y$-linear morphism (2.2.1) can be obtained by identifying $s$ with $-\dd_tt$ so that the action of $\dd_t$ is given by $-s\,t^{-1}$. The bijectivity of (2.2.1) can be shown by using increasing filtrations by the order of $s$ and $\dd_t$.
\sk
As in the proof of \cite[Proposition 2.5]{MP3}, we have the equalities for $i\in\N$
$$\dd_t^it^i=Q_i(\dd_tt)\q\h{with}\q Q_i(x):=x(x+1)\cdots(x+i-1).
\leqno(2.2.4)$$
(Consider, for instance, their actions on $t^{\al}$ to verify the equality.) These imply that
$$\dd_t^j\delta(t{-}f)=f^{-j}Q_j(-s)f^s.
\leqno(2.2.5)$$
So the $Q_i(x)$ give an explicit description of the isomorphism (2.2.1).
\msn
{\bf Remark~2.2}\,(i). The action of a vector field $\xi\in\Theta_X$ on $(j_f)_+M=M[\dd_t]\delta(t{-}f)$ is {\it twisted\1} as in (1.3.1), and this is compatible with (2.2.7). We can explain this twist by considering two lifts of vector fields $\xi$ to $Y=X\times\C$ as follows: The first one is the natural one using the product structure of $X\times\C$, and is denoted by $\xi^{(1)}$. The second one is defined by
$$\xi^{(2)}:=\xi^{(1)}+\xi(f)\dd_t.$$
We have
$$\xi^{(2)}f=0,\q\h{hence}\q\xi^{(2)}\delta(t{-}f)=0.$$
This implies the twist of the action of vector fields as in (1.3.1), since
$$\xi^{(1)}=\xi^{(2)}-\xi(f)\dd_t.$$
This argument can be extended to the case $M$ is any holonomic $\D$-module.
\sk
Here it may be better to use an {\it trivialization\1} of the projection $Y\to X$ given by the section defined by the graph of $f$ so that the direct image as $\D$-module is defined in the {\it product case}. (One can also use an automorphism of $Y$ over $X$ defined by $(x,t)\mapsto(x,t{-}f(x))$ and its inverse.)
\msn
{\bf Remark~2.2}\,(ii). In \cite[(1.1)]{rp}, two coordinate systems of $Y$ are used to explain the above twist of the action of vector fields: the first one is the natural one, that is, $(x_1,\dots,x_n,t)$, and the second one is given by $(x_1,\dots,x_n,t{-}f)$. (Note that the vector field $\dd_{x_i}$ {\it heavily depends on the choice of other coordinates}.)
\msn
{\bf 2.3.~The second isomorphism.} For $\al\in\Q$, it is also well-known to specialists that there are isomorphism of $\D_X\langle s,t\rangle$-modules (on which the action of $t$ is bijective)
$$\begin{array}{ccl}\OO_X(*Z)[s]f^s&\simto&\bl(\OO_X(*Z)f^{-\al}\br)[s]f^s\\
\rlap{$\hskip.26pt{\scriptstyle\,|}$}{\scriptstyle\cup}& &\q\q\q\rlap{$\hskip.26pt{\scriptstyle\,|}$}{\scriptstyle\cup}\\ \,\,g\1s\1^j\!f^s&\mapsto&\,\,\,\,g\1f^{-\al}\1(s{-}\al)^j\!f^s\,\,\,\,\bl(\,g\in\OO_X(*Z)\br).\end{array}
\leqno(2.3.1)$$
This is $\D_X\langle s,t\rangle$-{\it linear up to the shift of $s$ by $-\al$} (more precisely, it is $\D_X\langle s,t\rangle$-linear when the action of $\D_X\langle s,t\rangle$ on the target is given via an automorphism of $\D_X\langle s,t\rangle$ which is the identity on $\D_X\langle t\rangle$ and sends $s$ to $s{-}\al$. It should be noted that this shift of $s$ is required to get the {\it compatibility with the action of} $\D_X$ (or the vector field $\xi$) as is explained in (2.2.3). This implies that the $V$-filtration is {\it shifted by} $-\al$ under the isomorphism (2.3.1) using the {\it uniqueness\1} of $V$-filtration.
\sk
Note that the second isomorphism (2.3.1) for $\al=1$ coincides with the action of $t^{-1}$.
\msn
{\bf Remark~2.3.} Concerning the second isomorphism (2.3.1), some technical difficulty was noted in \cite[Section 1.8]{DS1}. However, this does not seem to cause any problems related to Hodge ideals for $\Q$-divisors (especially the determination of minimal exponents) as long as examples are computed explicitly.
\msn
{\bf 2.4.~Formula of Musta\c{t}\v{a} and Popa \cite[Theorem A]{MP3}.} In the notation of (2.1--3), the above two isomorphisms (2.2.1) and (2.3.1) immediately imply the following.
\msn
{\bf Theorem~2.4} (see \cite[Theorem A]{MP3} for the algebraic case). {\it Assume {\rm (2.1.7) (}that is, $\lceil D\rceil=Z)$. Then for $p\in\N$, we have the equalities
$$I''_p(D)=\bl\{\msum_{j=0}^p\,Q_j(\al)f^{p-j}g_j\mid\msum_{j=0}^p\,g_j\1\dd_t^j\1\delta(t{-}f)\in V^{\al}(i_f)_+\OO_X\br\},
\leqno(2.4.1)$$
where $g_j\in\OO_X$, and $I''_p(D)$ is as in {\rm (2.1.8--9)}.}
\msn
{\it Proof.} Setting
$$(M,F):=(\OO_X(*Z)f^{-\al},F),\q\q(M_f,F):=(i_f)_+(M,F),$$
we have the description of the Hodge filtration $F$ on $M_f$ as in \cite[(3.2.3.2)]{mhp}:
$$F_pM_f=\msum_{i\ges 0}\,\dd_t^i\bl(V^0M_f\cap j_*j^{-1}F_{p-i}M_f\br)\q(p\in\Z),
\leqno(2.4.2)$$
where $j:X\times\C^*\into X\times\C$ is the inclusion. (Note that $M_f$ does not mean the localization of $M$ by $f$, which coincides with $M$ itself by the bijectivity of the action of $f$ on $M$.)
\sk
By definition the direct sum decomposition
$$(i_f)_+M=\mopl_{p\in\N}\,M\,\dd_t^p
\leqno(2.4.3)$$
is compatible with the Hodge filtration $F$. (Note that the restriction of $F$ over $X\times\C^*$ is given by the order of $\dd_t$.) We then get
$$F_pM={\rm pr_0}(F_p(i_f)_+M)\q(p\in\Z),
\leqno(2.4.4)$$
using the projection associated with the direct sum decomposition (2.4.3)
$${\rm pr_0}:(i_f)_+M\onto M\,(=M\,\dd_t^0).
\leqno(2.4.5)$$
This implies that it is enough to consider the case $i=0$ in the right-hand side of (2.4.2).
\sk
Note that the projection ${\rm pr}_0$ is {\it compatible\1} with the first isomorphism (2.3.1), and can be given for $M[s]f^s$ by the {\it constant term\1} of a polynomial in $s$, that is, by {\it substituting\1} $0$ in $s$ (since $s=-\dd_tt$).
\sk
Returning to the proof of (2.4.1), the polynomial $Q_i(-s)$ in (2.2.5) is transformed into $Q_i(\al{-}s)$ by the second isomorphism (2.3.1), and we get $Q_i(\al)$ in (2.4.1) by substituting $s=0$ into $Q_i(\al{-}s)$. The equality (2.4.1) then follows from (2.4.4). Here we have $g_j\in\OO_X$, instead of $g_j\in\OO_X(*Z)$, since there are isomorphisms
$$V^{\al}(i_f)_+\OO_X\,\simto\,V^{\al}(i_f)_+(\OO_X(*Z))\q\q(\al>0).
\leqno(2.4.6)$$
This follows from \cite[Lemma 3.1.3 and Corollary 3.1.5]{mhp}, where $V_{\al}=V^{-\al}$. This finishes the proof of Theorem~(2.4).
\sk
As a corollary of Theorem~(2.4), we get the following.
\msn
{\bf Corollary~2.4} (see \cite[Theorem A]{MP3} for the algebraic case). {\it Under the assumption of Theorem~$(2.4)$, we have the equalities}
$$I''_p(D)=\V^{\al+p}\OO_X\,\,\,\,\,\h{mod}\,\,\,\,\,(f)\,\q(\forall\,p\in\N).
\leqno(2.4.7)$$
\msn
{\it Proof.} This follows from Theorem~(2.4) by the definition of microlocal $V$-filtration in (1.3.3), since $D=\al({\rm div}\1f)$, $\lceil D\rceil=Z$, hence $\al\in(0,1]$ (and $Q_i(\al)\ne 0$ for $\al>0$). This finishes the proof of Corollary~(2.4).
(This was shown in \cite{hi} when $\al=1$.)
\msn
{\bf Remark~2.4}\,(i). The above reasoning is expected to simplify some arguments related to $Q_i(x)$ in Appendix of the original version of \cite{MP3}.
\msn
{\bf Remark~2.4}\,(ii). Restricting to a neighborhood of a smooth point of $Z$, assume $f=x_1^m$ $(m\ges 2$) with $x_1$ a local coordinate of $X$. Then
$$\V^{\al'}\OO_X=x_1^k\OO_X\q\h{if}\q\lceil\al'm\rceil=k+1+\lfloor k/(m{-}1)\rfloor.
\leqno(2.4.8)$$
This follows from \cite[2.6]{MSS2} or \cite[2.2]{hi} using the compatibility of the microlocal $V\!$-filtration with smooth pull-backs.
Since ${\rm Sing}\,Z=\emptyset$ and ${\rm div}\1f=mZ$, we have by (2.1.8)
$$I''_p(\al mZ)=x_1^{p(m-1)}\OO_X\q\h{if}\q\al m\in(0,1].
\leqno(2.4.9)$$
These are compatible with (2.4.7). Indeed, setting $\al'=\al{+}p$, $k=p(m{-}1)$, we get
$$\lceil\al'm\rceil=pm+1=k+1+\lfloor k/(m{-}1)\rfloor.
\leqno(2.4.10)$$
\msn
{\bf Remark~2.4}\,(iii). In \cite[Theorems 0.4 and 0.6]{mos}, $\V^k\OO_X$, $\V^{>k}\OO_X$ for $k\ges 2$ should be replaced by the corresponding Hodge ideals. (See also \cite[4.1]{MSS2} for $V^{>1}\OO_X/\V^{>1}\OO_X$.)
\msn
{\bf 2.5.~Semi-weighted-homogeneous or non-degenerate case.} In the notation of (2.1), assume $D=\al Z$ ($\al\in(0,1]$), $Z={\rm div}\1f$ has an isolated singularity at 0, and moreover $f$ is semi-weighted-homogeneous with weights $w_i$, or has non-degenerate Newton boundary, for some local coordinates $x_1,\dots,x_n$ as in the introduction. In the non-degenerate case we assume $f$ is {\it convenient\1} (that is, the intersection of the Newton polyhedron $\Gp(f)$ with every coordinate axis of $\R^n$ is non-empty) so that \cite[Proposition B.1.2.3]{BGMM} can be applied. Set
$$\C\{x\}^{\ges\be}:=\msum_{v(g)\ges\be}\,\C\{x\}g\,\,\subset\,\,\C\{x\}=\OO_{X,0},
\leqno(2.5.1)$$
where $v(g)$ is as in (1.6--7). (It does not seem completely clear whether the two definitions of $v(g)$ in (1.6--7) coincide in the case $f$ is semi-weighted-homogeneous and has non-degenerate Newton boundary at the same time, although the right-hand side of (2.5.2) or (2.6,4) below should be the same.)
We have the following.
\msn
{\bf Theorem~2.5} (see \cite[Theorem 5.5]{Zh}). {\it Under the above assumptions, we have the following equalities for $\al\in(0,1]$, $p\in\N:$}
$$F_p\bl(\OO_{X,0}(*Z)f^{-\al}\br)=\msum_{j=0}^p\,F_{p-j}\D_{X,0}\bl(\C\{x\}^{\ges\al+j}f^{-\al-j}\br).
\leqno(2.5.2)$$
\sk
(This was shown in \cite[Theorem~0.9]{mos} in the semi-weighted-homogeneous case for $\al=1$.)
\msn
{\it Proof.} In the notation of (1.3), we can reduce (2.5.2) by the arguments in (2.2--4) to the following for $\al\in\Q$, $p\in\Z$:
$$F_pV^{\al}\bl(\OO_{X,0}[\dd_t,\dd_t^{-1}]\delta(t{-}f)\br)=\msum_{j\les p}\,F_{p-j}\D_{X,0}\bl(\C\{x\}^{\ges\al+j}\dd_t^{\1j}\delta(t{-}f)\br).
\leqno(2.5.3)$$
Indeed, the microlocalization morphism
$$\iota:\bl(\OO_X[\dd_t]\delta(t{-}f);F,V\br)\into\bl(\OO_X[\dd_t,\dd_t^{-1}]\delta(t{-}f);F,V\br)$$
induces the filtered isomorphisms for $\al<1:$
$$\Gr_V^{\al}\iota:\Gr_V^{\al}\bl(\OO_X[\dd_t]\delta(t{-}f),F\br)\simto\Gr_V^{\al}\bl(\OO_X[\dd_t,\dd_t^{-1}]\delta(t{-}f),F\br),$$
and we have for $\al\les 1$, $p\in\Z$
$$\iota^{-1}F_pV^{\al}\bl(\OO_X[\dd_t,\dd_t^{-1}]\delta(t{-}f)\br)=F_pV^{\al}\bl(\OO_X[\dd_t]\delta(t{-}f)\br).
\leqno(2.5.4)$$
(Here it is enough to show the assertion for $F$ and $V$ separately.) Note also that the projection ${\rm pr}_0$ in (2.4.5) is {\it compatible\1} with the action of $\D_X$, and $Q_i(\al)\ne 0$ for $\al>0$.
\sk
The assertion (2.5.3) is shown in \cite{mos} for the semi-weighted-homogeneous case and in \cite{Zh} for the non-degenerated case, where a non-trivial assertion \cite[Proposition B.1.2.3]{BGMM} is used and the convenience condition is needed, although it can be avoided later as is shown in Theorem~(2.7) below. (It is rather easy to show the inclusion $\supset$ in the assertions (2.5.2--3) without assuming the convenience condition in the non-degenerate case, see \cite[Proposition 3.2]{mic}.) This finishes the proof of Theorem~(2.5).
\msn
{\bf 2.6.~More explicit description of Hodge ideals.} To describe the Hodge ideals more explicitly, consider the differential operators
$$P(i,\be):=f\dd_{x_i}-\be f_i\,\in\,\D_{X,0}\q(i\in[1,n],\,\,\be\in\Q),
\leqno(2.6.1)$$
where $f_i:=\dd_{x_i}f$. This is justified by the equality
$$\dd_{x_i}(gf^{-\be})=\bl(P(i,\be)\1g\br)f^{-\be-1}.
\leqno(2.6.2).$$
For $\ib^{(k)}=(i^{(k)}_0,\dots,i^{(k)}_{k-1})\in[1,n]^k$ and $\al\in\Q$, set
$$\Pt(\ib^{(k)},\al):=P(i^{(k)}_{k-1},\al{+}k{-}1)\ssc\cdots\ssc P(i^{(k)}_1,\al{+}1)\ssc P(i^{(k)}_0,\al)\,\in\,\D_{X,0},$$
where $\Pt(\ib^{(k)},\al):=1$ for $k=0$. We have the equalities
$$\dd_{i_{k-1}^{(k)}}\ssc\cdots\ssc\dd_{i_0^{(k)}}(gf^{-\be})=\bl(\Pt(\ib^{(k)},\be)\1g\br)f^{-\be-k},
\leqno(2.6.3)$$
where $\dd_i:=\dd_{x_i}$. Theorem~(2.5) is then equivalent to the following formula for the Hodge ideals.
\msn
{\bf Corollary~2.6.} {\it Under the assumptions of Theorem~$(2.5)$, we have the following equalities for $\al\in(0,1]$, $p\in\N:$}
$$I_p(\al\1Z)=\msum_{k=0}^p\msum_{\ib^{(k)}\in[1,n]^k}\C\{x\}\1\Pt(\ib^{(k)},\al{+}p{-}k)\1\C\{x\}^{\ges\al+p-k}.
\leqno(2.6.4)$$
\sk
This is essentially equivalent to an inductive formula as in \cite[(10)]{Zh}. In the non-degenerate case, we can avoid the {\it convenience\1} assumption as is explained below.
\msn
{\bf 2.7.~Non-convenient non-degenerate case.} In this subsection we show the following.
\msn
{\bf Theorem~2.7.} {\it In the non-degenerate Newton boundary case, the assumption that $f$ is convenient is not necessary for Theorem~$(2.5)$ and Corollary~$(2.6)$ if $f$ is a polynomial with an isolated singularity at $0$ $($using the finite determinacy of holomorphic functions with isolated singularities$)$.}
\msn
{\it Proof.} Set
$$f^{\bf a}_c:=f+c\,\msum_{i=1}^n\,x_i^{a_i}.$$
For any $m\in\Z$, there are integers $a_i\ges m\,\,(i\in[1,n])$ such that $f^{\bf a}_c$ has non-degenerate Newton boundary for any $c\in\De$, see Proposition~(A.2) in Appendix below. We can define the Newton filtration $V_c$ on $\OO_{X,0}$ by the condition: $v_c(g)\ges\al$ associated with the Newton polyhedron $\Gp(f^{\bf a}_c)$ as in (1.7) for any $c\in\De$ (even if $f^{\bf a}_0=f$ is {\it not\1} convenient when $c=0$).
Since $\Gp(f^{\bf a}_c)\supset\Gp(f^{\bf a}_0)$, we have the inclusions
$$V_c^{\be}\OO_{X,0}\supset V_0^{\be}\OO_{X,0}=\C\{x\}^{\ges\be}\q(\be\in\Q).$$
\sk
Let $I^{\al+p}_c$ be the right-hand side of (2.6.4) associated with $f^{\bf a}_c$ ($c\in\De$), where $\C\{x\}^{\ges\al+p-k}$ is replaced by $V_c^{\ges\al+p-k}\OO_{X,0}$. If $m$ is sufficiently large ({\it depending on fixed\1} $\al$, $p$), we can prove the equalities\,:
$$I^{\al+p}_c=I^{\al+p}_0\q\h{in}\,\,\,\OO_{X,0}\q(\forall\,c\in\De,\,\al\in(0,1],\,p\in\N).
\leqno(2.7.1)$$
For the proof, we need the following.
\msn
{\bf Sublemma~2.7.} {\it An ideal $I\subset\OO_{X,0}$ contains $\mm_{X,0}^k$ if this holds mod $\,\mm_{X,0}^{k+1}$.}
\ms
(This follows from Nakayama's lemma applied to a subideal of $I$ whose image in $\OO_{X,0}/\mm_{X,0}^{k+1}$ coincides with $\mm_{X,0}^k/\mm_{X,0}^{k+1}$.)
\sk
Using Sublemma~(2.7), we may consider everything mod $\mm_{X,0}^{k+1}$ for $k$ very large, and we can take the above $m$ much larger than this $k$.
More precisely, since $f$ has an isolated singularity, there is an integer $k(\al,p)>\deg f$ such that $I_p(\al Z)\supset\mm_{X,0}^{k(\al,p)}$ and $f$ is $k(\al,p)$-determined, that is, $g\in\OO_{X,0}$ is {\it right equivalent\1} to $f$ if $f-g\in\mm_{X,0}^{k(\al,p)+1}$. These imply that $I_c^{\al+p}\supset\mm_{X,0}^{k(\al,p)}$ for $m>k(\al,p)$, $c\in\De^*$, using Corollary~(2.6). We can then prove (2.7.1) for $m\gg k(\al,p)$ with $\al,p$ fixed, since it is enough to show the equality $I_c^{\al+p}=I_0^{\al+p}$ mod $\mm_{X,0}^{k(\al,p)+1}$ by the Sublemma~(2.7). Note that for $j\in[0,p]\cap\Z$, $m\gg k(\al,p)$, we have
$$\bl\{\nu\in\N^n\mid|\nu|\les k(\al,p)+n\br\}\cap(\al{+}j)\bl(\Gp(f^{\bf a}_c)\setminus\Gp(f)\br)=\emptyset.
\leqno(2.7.2)$$
\sk
Consider now the one-parameter family $\{f^{\bf a}_c\}_{c\in\De}$. If $m\gg 0$, this should be an {\it analytically trivial\1} family (after an automorphism of $X\times\De$ over $\De$) by finite determinacy of holomorphic functions with isolated singularities and also by the existence of miniversal unfoldings.
Set
$$Z_c:=\{f^{\bf a}_c=0\}\subset X,\q \widetilde{Z}:=\h{$\bigsqcup$}_c\,Z_c\subset X\times\De.$$
By the analytic triviality of the family, we get the equalities
$$I_p(\al Z_c)=I_p(\al\widetilde{Z})|\!|_{u=c}\q\h{in}\,\,\,\OO_X\q(\forall\,c\in\De),
\leqno(2.7.3)$$
with $u$ the coordinate of $\De$.
The right-hand side denotes the image of $I_p(\al\widetilde{Z})$ in $\OO_X$ by identifying $\{u=c\}\subset X\times\De$ with $X$, and is denoted by $I_p(\al\widetilde{Z})\cdot\OO_{X\times\{c\}}$ in \cite{MP3}.
In our case this is essentially the same as the usual restriction as $\OO$-modules by applying the tensor product with $\OO_{\De}/(u{-}c)\OO_{\De}$ over $\OO_{\De}$ (or the snake lemma for the multiplication by $u{-}c$) to the short exact sequence
$$0\to I_p(\al\widetilde{Z})\to\OO_{X\times\De}\to\OO_{X\times\De}/I_p(\al\widetilde{Z})\to 0,$$
since the $\OO_{X\times\De}/I_p(\al\widetilde{Z})$ are free $\OO_{\De}$-modules of finite type by the analytic triviality.
\sk
We thus get a deformation of $\mm_{X,0}$-primary ideals of $\OO_{X,0}$ depending on the coordinate of $\De$ holomorphically, in particular, {\it continuously}. Here we may consider the finite dimensional vector spaces $I_p(\al Z_c)/\mm_{X,0}^{k(\al,p)}$, instead of $I_p(\al Z_c)\supset\mm_{X,0}^{k(\al,p)}$, together with the short exact sequences of finite dimensional vector spaces
$$0\to I_p(\al Z_c)/\mm_{X,0}^{k(\al,p)}\to\OO_{X,0}/\mm_{X,0}^{k(\al,p)} \to\OO_{X,0}/I_p(\al Z_c)\to 0,$$
as well as the corresponding exact sequences of finite free $\OO_{\De}$-modules.
\sk
By Corollary~(2.6), we have the equalities
$$I_p(\al Z_c)=I^{\al+p}_c\q\h{in}\,\,\,\OO_{X,0}\q(\forall\,c\in\De^*),
\leqno(2.7.4)$$
and these hold also for $c=0$ by (2.7.1), (2.7.3). This finishes the proof of Theorem~(2.7).
\msn
{\bf Remark~2.7}\,(i). The above argument cannot be extended to the proof of a generalization of \cite[Corollary~B]{Zh} for non-degenerate functions, since \cite{BGMM} would be needed again for the passage from Corollary~(2.6).
\msn
{\bf Remark~2.7}\,(ii). It is known that $f$ is $k$-determined if $k\ges\mu_f+1$, more precisely, if $\mm_{X,0}^{k+1}\subset\mm_{X,0}^2(\dd f)$, see for instance \cite[Theorem~9.1.4]{dJP}.
\msn
{\bf Remark~2.7}\,(iii). If $f$ has an isolated singularity at 0 and $f\in\mm_{X,0}^3$, then, for each $i\in[1,n]$, there is $j\in[1,n]$ such that
$$a_i\ee_i+\ee_j\in\Spx f\q\h{with}\q a_i\ges 2.$$
(In this paper $\ee_i\in\RR^n$ denotes the $i$\1th unit vector.) For instance, $f=\msum_{i=1}^n\,x_i^ax_{i+1}$ ($a\ges 2$) has an isolated singularity at 0, where $i$ is considered mod $n$ so that $x_{n+1}=x_1$.
\msn
{\bf Remark~2.7}\,(iv). The combinatorial data of the compact faces of $\Gp(f^{\bf a}_c)$ (with $c\in\C^*$) may heavily depend on the choice of the $a_i$. It does not seem clear whether the assertion at the beginning of the proof of Theorem~(2.7) holds with $a_i=m$ ($\forall\,i\in[1,n]$), since under the last hypothesis, any compact face $\si$ of the Newton polyhedron $\Gp(f^{\bf a}_c)$ is not necessarily the convex hull of the union of a compact face $\si'$ of $\Gp(f)$ and $m\1\ee_i$ ($i\in I$) for some $I\subset\{1,\dots,n\}$ with $|I|=\dim\si-\dim\si'$. For instance, if $f=\msum_{i=1}^{n-1}\,x_ix_{i+1}$ with $n$ even, then $\Gp(f)$ has a unique $(n{-}2)$-dimensional compact face $\si$, and $\Gp(f^{\bf a}_c)$ with $a_i=m$ ($\forall\,i\in[1,n]$) and $c\in\C^*$ has two $(n{-}1)$-dimensional compact faces which are respectively the convex hulls of $\si\cup\{\ee_{2j-1}\}_{1\les j\les n/2}$, $\si\cup\{\ee_{2j}\}_{1\les j\les n/2}$, and are contained in the hyperplanes defined by
$$\msum_{j=1}^{n/2}\bl(\nu_{2j-1}+(m{-}1)\nu_{2j}\br)=m,
\q\msum_{j=1}^{n/2}\bl((m{-}1)\nu_{2j-1}+\nu_{2j}\br)=m.$$
In the case $n=4$, the number of $k$-dimensional compact faces of $\Gp(f^{\bf a}_c)$ for $k=0,1,2,3$ with $a_i=m\ges 3$ and $c\in\C^*$ seems to be $7, 14, 11, 3$ respectively. The first two $3$-dimensional compact faces are both half of an octahedron, and intersect each other along a triangle (the intersection with some coordinate hyperplane is a rectangle). The third $3$-dimensional compact face is a tetrahedron in the hyperplane $\{\nu_1+(m{-}1)\nu_2+(m{-}1)\nu_3+\nu_4=m\}$.
\bs\bs
\vbox{\centerline{\bf 3. Proofs of the main theorems and propositions}
\bsn
In this section we prove the main theorems and propositions.}
\msn
{\bf 3.1.~Proof of Theorem~1.} It is enough to show the equality (10). Since the minimal or maximal exponent has multiplicity 1 (see (1.1.4)), we get by condition~(7)
$$\dim f\,\Om_f^n=\dim V^{\al_{f,\,\mu_f}}\,\Om_f^n=1.
\leqno(3.1.1)$$
(Recall that $\Om_f^n\cong\C\{x\}/(\dd f)$.) This implies that
$$\mm_{X,0}\1 f\,\Om_f^n=0.
\leqno(3.1.2)$$
\sk
Assume $\ep_f>0$. We first show the inequality
$$\al_{f,\,\mu_f}^{\HI}\ges\ga_f+1,$$
or equivalently, the inclusion
$$f\,\Om_f^n\subset V_{\HI}^{\ga_f+1}\Om_f^n.
\leqno(3.1.3)$$
By the definition of $\ga_f$ together with (2.4.7) and Proposition~(1.4), there is $g\in\mm_{X,0}\setminus\mm_{X,0}^2$ satisfying
$$g\in I_p(\al Z),\q\h{that is,}\q gf^{-\al-p}\in F_p(\OO_X(*Z)f^{-\al}),
\leqno(3.1.4)$$
where $\ga_f=\al{+}p$ with $\al\in(0,1],\,\,p\in\N$. (Note that $(\dd f)\subset\mm_{X,0}^2$.)
\sk
There is $i\in[1,n]$ such that $\dd_{x_i}g\in\C\{x\}$ is invertible.
Applying $\dd_{x_i}$ to $gf^{-\al-p}$, we then get that
$$f\in I_{p+1}(\al Z).$$
Indeed, the inclusion $\dd_{x_i}F_p\subset F_{p+1}$ together (2.6.1--2) implies that
$$\aligned&\bl(P(i,\al{+}p)\1g\br)f^{-\al-p-1}\in F_{p+1}(\OO_X(*Z)f^{-\al}),\\&\h{with}\q\q P(i,\al{+}p)\1g=(\dd_{x_i}g)f\,\,\,\,\h{mod}\,\,\,\,(\dd f).\endaligned$$
Thus (3.1.3) is proved. (Note that condition~(A) is not needed for the proof of (3.1.3).)
\sk
We now show the non-inclusion
$$f\,\Om_f^n\not\subset V_{\HI}^{\be}\Om_f^n\q\h{for}\q\be>\ga_f+1.
\leqno(3.1.5)$$
Since $f\in\mm_{X,0}^3$ and $\ga_f+1>\al_{f,\,\mu_f}$, this assertion can be reduced to the following.
$$\ga_f=\max\bl\{\be\in\Q\mid V^{\be}\Om_f^n\not\subset\mm_{X,0}^2\Om_f^n\br\},
\leqno(3.1.6)$$
using Corollary~(2.6) (where condition~(A) is needed). Indeed, we have by (3.1.2)
$$(\dd_ig)f\Om_f^n=0\q\h{if}\q g\in\mm_{X,0}^2,$$
and the assumption $f\in\mm_{X,0}^3$ implies that
$$P(i,\be)g=(\dd_ig)f-\be gf_i\in\mm_{X,0}^{k+2}\q\h{if}\q g\in\mm_{X,0}^k.
\leqno(3.1.7)$$
(If $f\notin\mm_{X,0}^3$, there is a problem in the $g$ constant case, see the proof of Proposition~1 in (3.4) below.)
\sk
We then get (3.1.5) and the equality (10) in the case $\ep_f>0$ (that is, $\al_{f,\,\mu_f}^{\HI}=\ga_f+1$), since (3.1.6) follows from the definition of $\ga_f$.
\sk
The argument is similar in the case $\ep_f\les 0$ (that is, $\ga_f+1\les\al_{f,\,\mu_f}$). Note first that the filtration $V$ on $\Om_f^n\cong\C\{x\}/(\dd f)$ is induced by the filtration $\C\{x\}^{\ges\be}$ ($\be\in\Q$) defined above. By Corollary~(2.6), we then get the inequality
$$\al_{f,\,\mu_f}^{\HI}\ges\al_{f,\,\mu_f},\q\h{or equivalently,}\q f\,\Om_f^n\subset V_{\HI}^{\al_{f,\,\mu_f}}\Om_f^n.$$
We can show the non-inclusion
$$f\,\Om_f^n\not\subset V_{\HI}^{\be}\Om_f^n\q\h{for}\q\be>\al_{f,\,\mu_f},$$
using Corollary~(2.6), since (3.1.6--7) hold also in the case $\ga_f+1\les\al_{f,\,\mu_f}$. So the equality (10) in the case $\ep_f\les 0$ (that is, $\al_{f,\,\mu_f}^{\HI}=\al_{f,\,\mu_f}$) follows. This finishes the proof of Theorem~1.
\msn
{\bf 3.2.~Proof of Theorem~2.} For the proof of (11) in Theorem~2, it is enough to show the inclusion (3.1.3), since $f\1\Om_f^n\subset\Om_f^n$ is a $\C\{x\}$-submodule generated by $[f\ddd x]$ and
$$\dim_{\C}f\1\Om_f^n=\mu_f-\tau_f=n-|I|.$$
So the argument is the same as in (3.1), and Theorem~2 follows. (Note that the arguments corresponding to (3.1.5--7) are not needed in the proofs of Theorems~2--3.) 
\msn
{\bf 3.3.~Proofs of Theorem~3.} We prove (3.1.3) with $f$, $\ga_f$ replaced by $fg$, $\ga_f(g)$ respectively, where we use (3.1.4) with $g$ replaced by $x_i\1g$. Here we need the condition that $g$ is a {\it monomial,} since this implies that $\dd_{x_i}(x_i\1g)=cg$ with $c\in\C^*$. This finishes the proof of Theorem~3.
\msn
{\bf 3.4.~Proof of Proposition~1.} We have $f\in\mm_{X,0}^2\setminus\mm_{X,0}^3$ by assumption. So there are $i,j\in[1,n]$ such that
$$\dd_{x_i}\dd_{x_j}f\,\,(=\dd_{x_i}f_j)\,\,\,\h{is invertible.}
\leqno(3.4.1)$$
It is well-known (see \cite{MP2}, \cite{MP3}) that
$$1\in I_p(\al Z)\q\h{with}\q\al{+}p=\al_{f,1}\,\,\,\bl(\al\in(0,1],\,p\in\N\br).
\leqno(3.4.2)$$
(This also follows from (2.4.7) and Proposition~(1.4).)
By (2.6.1--2) this implies that
$$f_j\in I_{p+1}(\al Z),\q\h{and then}\q f\1\dd_{x_i}f_j\in I_{p+2}(\al Z)\,\,\,\h{mod}\,\,\,(\dd f).
\leqno(3.4.3)$$
Here $f\notin(\dd f$), since $\mu_f\ne\tau_f$. So Proposition~1 follows.
\msn
{\bf 3.5.~Proof of Proposition~2.} The argument is similar to the proof of Proposition~1. Here $i=j=n$, and we have
$$g\in I_p(\al Z)\q\h{with}\q\al{+}p=v(g)\,\,\,\bl(\al\in(0,1],\,p\in\N\br).
\leqno(3.5.1)$$
using the easy part of Corollary~(2.6). We then get
$$g(2x_n)\in I_{p+1}(\al Z),\q\h{and}\q 2fg\in I_{p+2}(\al Z)\,\,\,\h{mod}\,\,\,(\dd f).
\leqno(3.5.2)$$
So Proposition~2 is proved.
\bs\bs
\vbox{\centerline{\bf 4. Examples}
\bsn
In this section we calculate some examples.}
\msn
{\bf 4.1.~Example.} Let
$$f=x^a+y^b+x^{a'}y^{b'},$$
with
$$a\ges b\ges 3,\q a'\in\bl[\tfrac{a-1}{2},a{-}2\br],\q b'\in\bl[\tfrac{b-1}{2},b{-}2\br],\q \tfrac{a'}{a}+\tfrac{b'}{b}>1.
\leqno(4.1.1)$$
This is semi-weighted-homogeneous with weights $\frac{1}{a}$, $\frac{1}{b}$. Here $\mu_f=(a-1)(b-1)$, and we have
$$\mu_f-\tau_f=(a-a'-1)(b-b'-1).
\leqno(4.1.2)$$
Indeed, we first see that the vector space $\C\{x,y\}/(\dd f)$ has a $\C$-basis consisting of 
$$[x^iy^j]\q\q(i,j)\in[0,a{-}2]{\times}[0,b{-}2],$$
using the $\mu$-constant deformation
$$f_u:=x^a+y^b+u\1x^{a'}y^{b'}\q(u\in\C),$$
together with the $\C^*$-action associated with the weights of $f$.
(Indeed, the deformation argument implies that the above assertion holds with $f$ replaced by $f_u$ for $|u|$ sufficiently small.)
\sk
Moreover, the $\C\{x,y\}$-submodule of $\C\{x,y\}/(\dd f)$ generated by $[f]=[c\1x^{a'}y^{b'}]$ ($c\in\C^*$) is annihilated by
$x^{a-1-a'}$, $y^{b-1-b'}$. This can be verified by using the assumption on $a',b'$ in (4.1.1).
So the equality (4.1.2) follows.
\sk
In this example, we see that
$$\al_{f,1}=\tfrac{1}{a}+\tfrac{1}{b},\q\ga_f=\tfrac{1}{a}+\tfrac{2}{b},\q\ep_f=\tfrac{2}{a}+\tfrac{3}{b}-1=\tfrac{6-(a-2)(b-3)}{ab}.
\leqno(4.1.3)$$
We thus get examples with $\ep_f$ positive or negative or 0 in Theorem~1, assuming
$$a'=a{-}2,\,\,\,b'=b{-}2\q\h{with}\q\tfrac{1}{a}+\tfrac{1}{b}<\tfrac{1}{2},$$
so that $\dim_{\C}f\1\Om_f^3=1$. (The last condition on $a,b$ is equivalent to that $\tfrac{a-2}{a}+\tfrac{b-2}{b}>1$.)
\msn
{\bf Remark~4.1.} Assume $(a,b)=(5,4)$ so that $\ep_f=\tfrac{3}{20}>0$. Then Theorem~1 implies that
$$\bl\{\exp\bl(-2\pi i\1\al_{f,j}^{\HI}\br)\,\big|\,j\in[1.\mu_f]\br\}\ne\bl\{\exp\bl(-2\pi i\1\al_{f,j}\br)\,\big|\,j\in[1.\mu_f]\br\}.
\leqno(4.1.4)$$
(This answers a question of the referee.) Note that the right-hand side coincides with the set of Milnor monodromy eigenvalues:
$$\bl\{\exp\bl(2\pi i(p/5\1{+}\1 q/4)\br)\,\big|\,p\in[1,4],q\in[1,3]\br\}.$$
\msn
{\bf 4.2.~Example.} Let
$$f=x^9+y^{10}+z^{11}+(x+y)x^3y^3z^3.$$
This is semi-weighted-homogeneous with weights $\frac{1}{9}$, $\frac{1}{10}$, $\frac{1}{11}$ (since $\frac{3}{9}+\frac{4}{10}+\frac{3}{11}=\tfrac{166}{165}>1$). In this example, we can show that the Hodge ideals $I_2(\al Z)$ are {\it not\1} weakly decreasing even modulo the Jacobian ideal $(\dd f)$ for $\al\in(0,1]$ as follows: We first see that
$$\aligned\dd_x\dd_x(x^2f^{-\al})&=\dd_x(2xf^{-\al}-\al x^2f_xf^{-\al-1})\\&=(2f^2-\al x^2f_{xx}f)f^{-\al-2}\mod \,\,(\dd f)f^{-\al-2},\endaligned
\leqno(4.2.1)$$
where $f_x:=\dd_xf$, $f_{xx}:=\dd_x^2f$, see also (2.6.3). Set
$$g_1=x^4y^3z^3,\q g_2=x^3y^4z^3,$$
so that $[f]$, $[x^2f_{xx}]$ are represented by $\C$-linear combinations of $[g_1],[g_2]$ in $\C\{x,y,z\}/(\dd f)$, more precisely,
$$[f]=-\tfrac{1}{990}\bl(17[g_1]+6[g_2]\br),\q[x^2f_{xx}]=-\bl(20[g_1]+18[g_2]\br).
\leqno(4.2.2)$$
\sk
Put
$$\aligned\be&:=v(x^2)=\tfrac{3}{9}+\tfrac{1}{10}+\tfrac{1}{11}=\tfrac{173}{330},\q \q\ep:=\tfrac{1}{90}.\\
\ga&:=v(g_1^2)=\tfrac{9}{9}+\tfrac{7}{10}+\tfrac{7}{11}=\tfrac{111}{330}+2\,\,\bl(>\,v(g_1g_2)\,>\,v(g_2^2)\br).\endaligned$$
Assume $\al\in(\be-\ep,\,\be]$. We have the linear independence of $[f^2]$, $[x^2f_{xx}f]$ in
$$\C\{x,y,z\}/\bl((\dd f)+\C\{x,y,z\}^{>\ga}\br),
\leqno(4.2.3)$$
since the $[x^iy^jz^k]$ for $(i,j,k)\in[0,7]{\times}[0,8]{\times}[0,9]$ form a $\C$-basis of $\C\{x,y,z\}/(\dd f)$ using the $\mu$-constant deformation
$$x^9+y^{10}+z^{11}+u\1x^4y^3z^3+v\1x^3y^4z^3\q(u,v\in\C),
\leqno(4.2.4)$$
together with the $\C^*$-action associated with the weights of $x,y,z$. (We can also see the linear independence of $[f^2]$, $[x^2f_{xx}f]$ in $\C\{x,y,z\}/(\dd f)$ using a computer.)
\sk
By Corollary~(2.6) it is then enough to show the following inclusions in the notation of (2.5):
$$\msum_i\,P(i,\al{+}1)\1\C\{x,y,z\}^{\ges\al+1}\,\subset\,\C\{x,y,z\}^{>\ga},
\leqno(4.2.5)$$
$$\msum_{i,j}\,P(i,\al{+}1)P(j,\al)\1\mm_{X,0}^3\,\subset\,\C\{x,y,z\}^{>\ga},
\leqno(4.2.6)$$
since $\al>\be-\ep>\ga-2$ so that $\C\{x,y,z\}^{\ges\al+2}\subset\C\{x,y,z\}^{>\ga}$. Here any monomials of degree~2 {\it except for\1} $x^2$ are {\it not\1} contained in $\C\{x,y,z\}^{\ges\al}$ (since $\ep={\rm wt}\,x-{\rm wt}\,y$).
So Corollary~(2.6) and (4.2.5--6) imply that $I_p(\al Z)$ mod $(\dd f)+\C\{x,y,z\}^{>\ga}$ is spanned by $[2f^2-\al x^2f_{xx}f]$.
\sk
The inclusion (4.2.5) holds, since
$$\msum_i\,P(i,\al{+}1)\1\C\{x,y,z\}^{\ges\al+1}\,\subset\,\C\{x,y,z\}^{\ges\al+2-\frac{1}{9}},
\leqno(4.2.7)$$
(where $\frac{1}{9}$ is the maximum of the weights of $x,y,z$) and
$$\bl(\be-\ep+2-\tfrac{1}{9}\br)-\ga=\tfrac{173}{330}-\tfrac{1}{90}-\tfrac{1}{9}-\tfrac{111}{330}=\tfrac{13}{198}>0.
\leqno(4.2.8)$$
Using the inclusion $\mm_{X,0}^3\subset\C\{x,y,z\}^{\ges v(z^3)}$, the assertion (4.2.6) follows from the inequality 
$$\bl(v(z^3)+2-\tfrac{2}{9}\br)-\ga=\tfrac{1}{9}+\tfrac{1}{10}+\tfrac{4}{11}-\tfrac{2}{9}-\tfrac{111}{330}=\tfrac{8}{495}>0.
\leqno(4.2.9)$$
So the Hodge ideals are not weakly decreasing even modulo the Jacobian ideal.
\msn
{\bf Remark~4.2}\,(i). The above argument implies that the $I_2(\al Z)$ mod $(\dd f)$ are not weakly decreasing {\it even restricted to $\al\in(0,1)\cap\tfrac{1}{m}\1\Z$ with\1} $m=990$ (the order of the monodromy).
\msn
{\bf Remark~4.2}\,(ii). Using Corollary~(2.6), it is easy to see that the Hodge ideals $I_1(\al Z)$ ($\al\in(0,1]$) are {\it not weakly decreasing without taking\1} mod $(\dd f)$, for instance, if $f=x^a+y^b\,$ with $\,a>b\ges 3$. Indeed,
$$f-\al\1xf_x\,\notin\,I_1(\al Z)\,\ni\,f-\al\1yf_y\q\h{if}\q\tfrac{2}{a}{+}\tfrac{1}{b}<\al\les\tfrac{1}{a}{+}\tfrac{2}{b}\,(<1),$$
where $v(f)=v(xf_x)=v(yf_y)=\tfrac{1}{a}{+}\tfrac{1}{b}{+}1<\al{+}1$, see also \cite[Example 4.6]{Zh} for $(a,b)=(5,2)$ with $\al\in\bl(\tfrac{9}{10},1\br]$.
\msn
{\bf 4.3.~Example.} Assume
$$\al_{f,\,\mu_f}-\al_{f,1}=1,\q\mu_f\ne\tau_f.
\leqno(4.3.1)$$
Then $\al_{f,1}=\tfrac{n-1}{2}$, and $\tau_f=\mu_f-1$ by (1.1.4), (1.2.6). From condition~(4.3.1) we can deduce the {\it unimodality\1} of $f$ using the local injectivity of period maps via Brieskorn lattices at smooth points of $\mu$-constant strata (see \cite[Theorem~3.2]{per}) together with a description of the Brieskorn lattice $H''_f$ using a good section or an opposite filtration, see for instance \cite[Proposition~3.4]{bl}, \cite[Corollary~1.4]{dabl}. Indeed, (4.3.1) implies that there are free generators $\om_i$ ($i\in[1,\mu_f]$) of $H''_f$ over $\C\{\!\{\dd_t^{-1}\}\!\}$ such that $\om_i=u_i$ ($i\ne 1$) and $\om_1=u_1+\xi_f\1\dd_tu_{\mu_f}$ for $\xi_f\in\C$ (choosing an opposite filtration). Here the $u_i$ are free generators of $H''_f$ corresponding to a basis of the Milnor cohomology compatible with a fixed opposite filtration and satisfying $\dd_ttu_i=\al_{f,i}u_i$, see {\it loc.\,cit.} (Note that $\xi_f$ depends on the choice of $u_1,u_{\mu_f}$.) If there is a $\mu$-constant one-parameter family of holomorphic functions with isolated singularities satisfying (4.3.1), the $u_i$ can be chosen to be stable by parallel translation, and the complex number $\xi_f$ varies holomorphically.
\sk
By the classification of holomorphic functions with isolated singularities and {\it low modalities\1} (see for instance \cite{AGLV}), we then get the following {\it normal form\1} after an appropriate analytic coordinate change:
$$f=\begin{cases}x_1^{p_1}{+}\1 x_2^{p_2}{+}\1 c\1x_1^2x_2^2\,\,\,\,\bl(\tfrac{1}{\,\raise1.8pt\h{$\scriptstyle p_1$}}{+}\tfrac{1}{\,\raise1.8pt\h{$\scriptstyle p_2$}}{<}\tfrac{1}{2}\br)&\h{if}\,\,n\,{=}\,2,\\ x_1^{q_1}{+}\1 x_2^{q_2}{+}\1 x_3^{q_3}{+}\1 c\1x_1x_2x_3\1{+}\1\msum_{i=4}^n\,x_i^2\,\,\,\,\bl(\tfrac{1}{\,\raise1.8pt\h{$\scriptstyle q_1$}}{+}\tfrac{1}{\,\raise1.8pt\h{$\scriptstyle q_2$}}{+}\tfrac{1}{\,\raise1.8pt\h{$\scriptstyle q_3$}}{<}1\br)&\h{if}\,\,n\,{\ges}\,3\raise14pt\h{}.\end{cases}
\leqno(4.3.2)$$
Here $c\in\C^*$, and we assume
$$p_1\ges p_2\,\,(\ges 3),\q q_1\ges q_2\ges q_3\,\,(\ges 2).$$
Note that the sum of the first polynomial in (4.3.2) with $\msum_{i=3}^n\,x_i^2$ for $n\,{\ges}\,3$ is right equivalent to the second one with $(q_1,q_2,q_3)=(p_1,p_2,2)$ replacing $x_3$ with $x_3+\sqrt{-c}\,x_1x_2$, where $c$ is also modified. (The condition $\al_{f,\,\mu_f}-\al_{f,1}<1$ is equivalent to that $Z$ has an isolated singularity of type $A,D,E$ by a similar argument.)
\sk
For $f$ as in (4.3.2), we have by \cite{St2} (for $n\,{=}\,2$) and \cite{exp} (in general)
$$\Sp_f(t)=\begin{cases}t^{1/2}\,\bl(1\1 {+}\1 t^{1/2}{+}\1 t\1 {+}\msum_{i=1}^2\,\msum_{0<j<p_i}\,t^{\1j/p_i}\br)&\h{if}\,\,\,n\,{=}\,2,\\ t^{(n-1)/2}\,\bl(1\1 {+}\1 t\1 {+}\msum_{i=1}^3\,\msum_{0<j<q_i}\,t^{\1j/q_i}\br)&\h{if}\,\,\,n\,{\ges}\,3\raise14pt\h{}.\end{cases}
\leqno(4.3.3)$$
(This also follows from the information of the Milnor monodromy using condition (4.3.1) together with the symmetry of spectral numbers, and can be reduced to the case $n=2$ or 3 by the Thom-Sebastiani theorem.)
By Theorem~1 and Proposition~1, we get the inequality
$$\al_{f,\,\mu_f}^{\HI}-\al_{f,\,\mu_f}\ges\tfrac{1}{\,\raise1.8pt\h{$\scriptstyle p_2$}}\,\,\,\h{or}\,\,\,\,\tfrac{1}{\,\raise1.8pt\h{$\scriptstyle q_n$}}.
\leqno(4.3.4)$$
Here the equality holds in the case $n=2$ or $n=3$ with $q_3\ges 3$ by Theorem~1, and we may set $q_i=1$ for $i>3$ with $n>3$ by Proposition~1.
\msn
{\bf Remark~4.3}\,(i). The singularity defined by the polynomial in (4.3.2) for $n\ges 3$ is called $T_{q_1,q_2,q_3}$. The polynomials in (4.3.2) are the only {\it unimodal\1} singularities {\it satisfying\1} (4.3.1), see for instance \cite{AGLV}.
Even if condition~(4.3.1) can imply that $f$ is $T_{q_1,q_2,q_3}$ when $n\ges 3$, it is not necessarily easy to determine $q_1,q_2,q_3$ explicitly.
\sk
For instance, take any $g(x,y)\in(x,y)^a\setminus(x,y)^{a+1}\subset\C\{x,y\}$ ($a\ges 3$), and put
$$f=g(x,y)+z^b+(x^2-y^2)z\,\in\,\C\{x,y,z\}\q\bl(\tfrac{2}{a}+\tfrac{1}{b}<1\br).
\leqno(4.3.5)$$
Substituting $x=u{+}v$, $y=u{-}v$, and setting $h(u,v)=g(u{+}v,u{-}v)$, we get
$$f=h(u,v)+z^b+4uvz\,\in\,\C\{u,v,z\}.
\leqno(4.3.6)$$
Unless $g$ is divisible by $u^2$ or $v^2$, it seems that $f$ has an isolated singularity (as far as calculated). In this case $f$ should be $T_{q_1,q_2,q_3}$, although it is not necessarily easy to determine $q_1,q_2,q_3$ in general, see Remarks~(4.3)\,(ii) and (iii) below.
\msn
{\bf Remark~4.3}\,(ii) Assume $g(x,y)=x^2y^2$ with $b=5$, for simplicity, in the notation of Remark~(4.3)\,(i) above. We can see that $(q_1,q_2,q_3)=(5,4,4)$ by showing that the zero locus of $f$ is isomorphic to that of $T_{5,4,4}$, calculating the multiplicative structure of $\OO_{X,0}/((\dd f),f)$ (see \cite{MY}), although the non-degeneracy condition is {\it not\1} satisfied for the convex hull of $(4,0,0)$, $(0,4,0)$, where $f$ is as in (4.3.6) with $h=(u^2{-}v^2)^2$. Indeed, using for instance a computer program Singular \cite{DGPS}, we can get that
$$\aligned&\dim_{\C}\OO_{X,0}/\bl((\dd f),f\br)=11,\\&\dim_{\C}\OO_{X,0}/\bl((\dd f),f,m\br)=10\q\h{if}\q m=uv,uz,vz.\endaligned$$
For instance, the last equality with $m=uv$ can be obtained by typing as follow:
\ms
\vbox{\small\sf\verb#ring R=0,(u,v,z),ds; ideal J; poly f=(u^2-v^2)^2+z^5+4*u*v*z;#
\sk
\verb#J=(jacob(f),f,uv); vdim(groebner(J));#}
\msn
Here {\small\sf\verb#ds#} means that the calculation is done in the localization of the polynomial ring at $0$.
\sk
Set $A(f):=\OO_{X,0}/\bl((\dd f),f\br)$. The above calculation implies that $\dim_{\C}\OO_{X,0}[m]=1$, that is, $\OO_{X,0}[m]=\C[m]\ne 0$ in $A(f)$, for $m=uv,uz,vz$, where $[m]$ is the class of $m$ in $A(f)$. This means that the annihilator of $[m]$ is the maximal ideal $\mm_{X,0}\subset\OO_{X,0}$, and we get for instance $[u^2v]=[uv^2]=0$. A similar assertion holds for $m=u^3,v^3,z^4$, since $4[uv]=-5[z^4]$, etc. We can then deduce that there is a $\C$-basis of $A(f)$ consisting of
$$[u^i]\,\,\,(i\in[1,3]),\,\,\,[v^j]\,\,\,(j\in[1,3]),\,\,\,[z^k]\,\,\,(k\in[0,4]),$$
giving an isomorphism $A(f)\cong A(f')$ as $\C$-algebras by repeating the same calculation with $f$ replaced by $f':=u^4+v^4+z^5+4uvz=f+2u^2v^2$.
(This example shows that the restriction of $f$ to $x_n=0$ may have {\it non-isolated singularities\1} even if $f$ has non-degenerate Newton boundary unless $f$ is convenient.)
\msn
{\bf Remark~4.3}\,(iii). In the notation of Remark~(4.3)\,(i) above, assume $g=y^3$ or equivalently $h=(u{-}v)^3$. We see that $\mu_f=b{+}5$ at least for $b\les 100$ by the same argument as in Remark~(4.3)\,(ii) above. It seems that $(q_1,q_2,q_3)=(b,3,3)$, since we have $\dim_{\C}\OO_{X,0}[m]=1$ for $m=u^2,v^2,z^{b-1}$ ($b\les 100$) according to Singular. Indeed, the last calculation implies that $\dim_{\C}\Gr_G^kA(f)=\{i\in[1,3]\mid k<q_i\}$ for $k\ges 1$ with $G^kA(f):=\mm_{X,0}^kA(f)$ ($k\ges 0$).
\sk
If $g=(x{+}y)y^2$, that is, $h=2u(u{-}v)^2$, then we can verify that $\mu_f=2b{+}2$ at least for $b\les 100$, and it seems that $(q_1,q_2,q_3)=(b,b,3)$, since $\dim_{\C}\OO_{X,0}[m]=1$ for $m=u^2,v^{b-1},z^{b-1}$ ($b\les 100$).
This is rather surprising. (The details are left to the reader.)
Here one can also apply {\small\sf\verb#spprint(spectrum(f));#} to $f$ in (4.3.6) with $h=2u(u{-}v)^2$ and also with $h=u^3+v^b$ for the comparison of spectra (although the latter result also follows from (4.3.3)) after adding {\small\sf\verb#LIB "gmssing.lib";#} before {\small\sf\verb#Ring R=0#} in the above code (where the spectral numbers are shifted by $-1$ as in \cite{St3}.)
Note that $q_1,q_2,q_3$ are uniquely determined by the spectrum of $f$ (and the latter depends only on the zero-locus of $f$).
\bs\bs
\vbox{\centerline{\bf Appendix. Key to the proof of Theorem~(2.7)}
\bsn
In this Appendix, we prove the key Proposition~(A.2) to the proof of Theorem~(2.7) after recalling some basics of Newton polyhedra.}
\msn
{\bf A.1.~Equations defining Newton polyhedra.}
In the notation of (1.7), the Newton polyhedra $\Ga(f)$, $\Gp(f)$ of a polynomial $f$ are defined by a finite number of inequalities $\ell_k\ges 0$ with $\ell_k$ linear functions with constant terms. Taking a minimal set of linear functions, which is denoted by $\LF(f)$ or $\LFp(f)$, we have a one-to-one correspondence between this set and the $(n{-}1)$-dimensional faces of $\Ga(f)$ or $\Gp(f)$, where each $\ell_k\in\LF(f)$ or $\LFp(f)$ is unique up to multiplication by a positive number. It is rather easy to determine explicitly $\LF(f)$ for a polynomial $f$, since $\Spx f$ is a {\it finite set}.
\sk
We say that a linear function with a constant term $\ell(\nu)=\msum_{i=1}^n\,c_i\nu_i+c_0$ is {\it strictly positive\1} if $c_i>0$ ($\forall\,i\in[1,n]$), and {\it weakly positive\1} if $c_i\ges 0$ ($\forall\,i\in[1,n]$).
Let $\LF^{\rm sp}(f)$, $\LF^{\rm wp}(f)$ respectively denote the subset of $\LF(f)$ consisting of strictly positive and weakly positive linear functions with constant terms, and similarly for $\LFp^{\rm sp}(f)$, $\LFp^{\rm wp}(f)$.
\sk
For ${\bf b}:=(b_i)\in\Z_{>0}^n$ and $\ga\in\C^*$ sufficiently general, set
$$h^{\bf b}:=\mprod_{i=1}^n(1+\ga\1x_i^{b_i}).$$
(Here we can take ${\bf b}\in\Q_{>0}^n$ or even ${\bf b}\in\R_{>0}^n$, using (A.1.3) below.)
\msn
{\bf Lemma~A.1.} {\it We have the equalities}
$$\LF^{\rm sp}(f)=\LF^{\rm sp}(h^{\bf b}f)=\LFp^{\rm sp}(f),
\leqno{\rm(A.1.1)}$$
\vskip-6mm
$$\LF^{\rm wp}(h^{\bf b}f)=\LFp(f)=\LFp^{\rm wp}(f).
\leqno{\rm(A.1.2)}$$
\msn
{\it Proof.} The last equality of (A.1.2) follows from the stability of $\Gp(f)$ by the action of $u\in\R^n_{\ges 0}$ via addition.
For the proof of other equalities, note first that
$$\Spx h^{\bf b}f=\mcup_{J\subset\{1,\dots,n\}}\bl(\Spx f+\msum_{j\in J}\,b_j\1\ee_j\br),
\leqno{\rm(A.1.3)}$$
since $\ga$ is sufficiently general. A linear function $\ell$ with a constant term belongs $\LF(f)$ (up to multiplication by a positive number) if and only if the following two conditions hold:
$$\dim\,\langle\ell^{-1}(0)\cap\Spx f\rangle^{\rm aff}=n{-}1,\q\Spx f\subset\{\ell\ges0\}.$$
Here $\langle S\rangle^{\rm aff}$ for a subset $S\subset\R^n$ denotes the {\it affine subspace spanned by\1} $S$, that is, the smallest affine subspace of $\R^n$ containing $S$. (Choosing a point $p\in S$, this coincides with $p+V_{S,p}$, where $V_{S,p}\subset\R^n$ is the vector subspace spanned by $q-p$ for $q\in S$.)
A similar assertion holds for $\LFp(f)$ with $\Spx f$ replaced by $\Spx f+\R_{\ges 0}^n$ (although the latter is {\it not finite}).
\sk
These imply the first equality of (A.1.1), that is, $\LF^{\rm sp}(f)=\LF^{\rm sp}(h^{\bf b}f)$. Indeed, if there is $\ell\in\LF(h^{\bf b}f)$ such that
$$\aligned&\ell^{-1}(0)\cap\Spx h^{\bf h}f\,\not\subset\,\Spx f,\q\q\q\q\q\q\h{that is,}\\ &\ell^{-1}(0)\cap\bl(\Spx f+\msum_{j\in J}\,b_j\1\ee_j\br)\ne \emptyset\q\h{for some}\,\,\,J\ne\emptyset,\endaligned$$
then $\ell$ cannot be strictly positive, since we get $\nu\in\Spx h^{\bf b}f$ with
$$\ell(\nu)=0,\q\ell(\nu-b_j\ee_j)\ges 0\q\h{for each $j\in J$}.$$
Similarly we have $\LF^{\rm sp}(f)=\LFp^{\rm sp}(f)$ using $\Spx f+\R_{\ges 0}^n$ instead of $\Spx h^{\bf b}f$.
\sk
For weakly positive linear equations $\ell$ which are not strictly positive, we can apply an inductive argument on $n$, using the projections $\pi_i:\R^n\to\R^{n-1}$ defined by omitting the $i\1$th coordinate if the coefficient $c_i$ of $\ell$ vanishes $(i\in[1,n])$. Indeed, we have
$$\pi_i\bl(\Ga(f)\br)=\Ga(f_{(i)}),\q\pi_i\bl(\Gp(f)\br)=\Gp(f_{(i)})\q\h{in}\q\R^{n-1},
\leqno{\rm(A.1.4)}$$
where $f_{(i)}$ is defined by substituting $x_i=\ga_i$ into $f$ for a sufficiently general $\ga_i\in\De_{\ep}^*$ (with $\De_{\ep}^*$ a punctured disk of radius $\ep\ll 1$). Using the last equality of (A.1.2) and applying the first equality of (A.1.4) to $h^{\bf b}f$, we then get
$$\aligned\LFp(f)&\,=\,\LFp^{\rm sp}(f)\,\sqcup\,\mcup_{i=1}^n\,\pi_i^*\LFp(f_{(i)}),\\ \LF^{\rm wp}(h^{\rm b}f)&\,=\,\LF^{\rm sp}(h^{\rm b}f)\,\sqcup\,\mcup_{i=1}^n\,\pi_i^*\LF^{\rm wp}\bl((h^{\rm b}f)_{(i)}\br).\endaligned
\leqno{\rm(A.1.5)}$$
So the first equality of (A.1.2) follows by induction on $n$, using the last one of (A.1.1). This finishes the proof of Lemma~(A.3).
\msn
{\bf Remark~A.1}\,(i). Using the above one-to-one correspondence between linear functions with constant terms and $(n{-}1)$-dimensional faces, we see that any vertex of $\Gp(f)$ belongs to $\Spx f$. Hence any compact face of $\Gp(f)$ is a convex hull of a finite subset of $\Spx f$ in the notation of (1.7) by Lemma~(A.1).
\msn
{\bf Remark~A.1}\,(ii). The affine subspace $A_{\si}$ spanned by a {\it compact\1} face $\si\subset\Gp(f)$ cannot be stable by the action of $\ee_i$ via addition for any $i\in[1,n]$. Indeed, we have
$$\si=\mcap_{k\in J}\,\ell_k^{-1}(0)\cap\Gp(f),$$
where the $\ell_k$ ($k\in J$) are weakly positive linear functions with constant terms defining the $(n{-}1)$-dimensional faces of $\Gp(f)$ containing $\si$. If $A_{\si}$ is stable by the action of $\ee_i$, then so are the $\ell_k^{-1}(0)$ and also $\si$ by the above formula. Hence $\si$ cannot be compact.
\msn
{\bf A.2.~Key Proposition to the proof of Theorem~(2.7).} The following is needed for the proof of Theorem~(2.7).
\msn
{\bf Proposition~A.2.} {\it Let $f$ be a polynomial having non-degenerate Newton boundary at the origin. Then, for any $m\gg 0$ and $i\in[1,n]$, there is $a_i\ges m$ such that $g:=f+c\1x_i^{a_i}$ has non-degenerate Newton boundary at the origin for any $c\in\C^*$.}
\ms
For the proof, we may assume $n\ges 3$. We first show the following lemma (which does not seem completely trivial to non-specialists):
\msn
{\bf Lemma~A.2.} {\it Set $p_i:=a_i\ee_i\in\R^n$. Any compact face $\si$ of $\Gp(g)$ is either a compact face of $\Gp(f)$ $($if $p_i\notin\si)$ or the convex hull of $\tau\cup\{p_i\}$ with $\tau$ a compact face of $\Gp(f)$ such that $\ee_i$ is not contained in the vector subspace $V_{\tau}\subset\R^n$ spanned by $\tau$ $($if $p_i\in\si)$.}
\msn
{\it Proof.} Let $\si$ be a compact face of $\Gp(g)$. As is noted in Remark~(A.1) (i), this is a convex hull of a finite subset
$$\Xi\,\subset\,\Spx g=\Spx f\sqcup\{p_i\}.$$
\msn
{\bf Case 1}: $p_i\notin\Xi$, that is, $p_i\notin\si$ (since $a_i\gg\deg f$). We first show that there is an $(n{-}1)$-dimensional face
$$\si'\subset\dd\Gp(g)\q\h{satisfying}\q p_i\notin\si'\supset\si.
\leqno{\rm(A.5)}$$
In the case $\dim\si=n{-}1$, we can take $\si'=\si$. We may thus assume $\dim\si\les n{-}2$. Let $\eta'$ be an $(n{-}1)$-dimensional face of $\Gp(g)$ containing $\si$. If $p_i\notin\eta'$, we get (A.5). So we may assume $p_i\in\eta'$.
There is an $(n{-}2)$-dimensional face
$$\eta\subset\eta'\q\h{satisfying}\q p_i\notin\eta\supset\si.$$
(We have $\eta=\si$ if $\dim\si=n{-}2$.) Let $q$ be an interior point of $\eta$. (Note that $n\ges 3$.) Let $V'$ be the 2-dimensional vector subspace of $\R^n$ spanned by $q$ and $\ee_i$ (or $p_i$). Consider a real $1$-dimensional subset
$$\Ga_{V'}:=\Gp(g)\cap V'\,\subset\,V'\cong\R^2.$$
One connected component $\Ga_1$ of $\dd\Ga_{V'}\setminus\{q\}$ is the union of $p_i+\R_{\ges 0}$ and the interior of the segment $[p_i,q]\subset\R^n$ connecting $p_i$ and $q$. (Note that this segment is contained in $\eta'$.) Then $\si'$ is obtained as an $(n{-}1)$-dimensional face of $\Gp(g)$ containing an non-empty neighborhood of $q$ in $\dd\Ga_{V'}\setminus\Ga_1$. So (A.5) follows.
\sk
Consider a weakly positive linear function with a constant term
$$\ell'=\msum_{j=1}^n\,c_j\nu_j+c_0$$
corresponding to $\si'\subset\Gp(g)$ so that
$$\si'=\ell'{}^{-1}(0)\cap\Gp(g),\q\Gp(g)\subset\{\ell'\ges 0\}.$$
Since $p_i\in\Gp(g)\setminus\si'$, we have $\ell'(p_i)>0$. Set
$$\si'':=\si'\cap\Gp(f)=\ell'{}^{-1}(0)\cap\Gp(f).$$
Since $\ell'$ is weakly positive and $\ell'(p_i)>0$, we get $\ell'(p_i{+}\R_{\ges 0}^n)\subset\R_{>0}$. This implies that
$$\si'=\si'',$$
since $\Gp(g)$ is the convex hull of the union of $\Gp(f)$ and $p_i{+}\R_{\ges 0}^n$. Thus $\si'$ is a face of $\Gp(f)$, and so is $\si$ which is a face of $\si'$.
\msn
{\bf Case 2}: $p_i\in\Xi$. We have
$$\Xi=\Xi'\sqcup\{p_i\}\q\h{with}\q\Xi'\subset\Spx f.$$
Let $\tau\subset\si$ be the convex hull of $\Xi'$.
The affine subspace $A_{\tau}\subset\R^n$ spanned by $\tau$ (or $\Xi'$) does not contain $p_i$ if $a_i\gg 0$, since $\Xi'$ is a subset of a {\it finite\1} set $\Spx f$.
Moreover, $A_{\tau}\subset\R^n$ is {\it unstable\1} by the action of $\ee_i$ via addition, since so is $A_{\si}\supset A_{\tau}$, see Remark~(A.1) (ii). This implies that $p_i\notin A_{\tau}$ (since $a_i\gg 0$), hence $\tau=A_{\tau}\cap\si$ is a face of $\si$.
Then by the same argument as in Case 1, we see that $\tau$ is a face of $\Gp(f)$.
\sk
It now remains to show that $\ee_i\notin V_{\tau}$. Assume $\ee_i\in V_{\tau}$. The affine subspace $A_{\tau}\subset\R^n$ spanned by $\tau$ must then contain $\alpha\1\ee_i$ for some $\alpha\in\R\setminus\{0\}$, since $A_{\tau}$ is not stable by the action of $\ee_i$ as is shown above.
(The assertion is reduced to the case $V_{\tau}=\R^{r+1}$ with $A_{\tau}=\{\nu_{r+1}=1\}\subset\R^{r+1}$ after some coordinate change.)
So the affine subspace $A_{\pi_i(\tau)}\subset\R^{n-1}$ spanned by $\pi_i(\tau)$ must contain $0$, where $\pi_i:\R^n\to\R^{n-1}$ is the projection in (A.1.4). Moreover, $\pi_i(\tau)$ must be contained in the boundary of $\pi_i\bl(\Gp(f)\br)=\Gp(f_{(i)})\subset\R^{n-1}$ (since $\si\subset\dd\Gp(g)$ and $a_i\gg 0)$, where $f_{(i)}$ is as in (A.1.4). 
These contradict the weak positivity of the coefficients of the linear functions defining $\Gp(f_{(i)})\subset\R^{n-1}$ (considering an $(n{-}2)$-dimensional face containing $\pi_i(\tau)$). We thus get that $\ee_i\notin V_{\tau}$. This finishes the proof of Lemma~(A.2).
\msn
{\bf A.3.~Proof of Proposition~(A.2).} We have to show that, for any compact face $\si$ of $\Gp(g)$, $g_{\si}$ is non-degenerate, that is, (1.7.2) is satisfied. In the first case of Lemma~(A.2), this follows from the latter. 
In the second case, we have $\nu^{(k)}\in\Z^n$ ($k\in[1,n]$) such that $\sum_{k=1}^n\Z\1\nu^{(k)}\subset\Z^n$ has a finite quotient group (that is, $\det(\nu^{(1)},\dots,\nu^{(n)})\ne 0$), $\nu^{(1)},\dots,\nu^{(r)}$ are free generators of $V_{\tau,\Z}:=V_{\tau}\cap\Z^n$ with $r:=\dim V_{\tau}$, and $\nu^{(n)}=\ee_i$. These define a finite \'etale morphism of affine tori
$$\rho:(\C^*)^n\to(\C^*)^n,$$
such that
$$\rho^*y_k=x^{\nu^{(k)}}:=\mprod_{j=1}^n\,x_j^{\nu_j^{(k)}},\q\h{that is,}\q\rho^*\log y_k=\msum_{j=1}^n\,\nu_j^{(k)}\log x_j,$$
where the $x_j$ and $y_k$ are the coordinates of the tori. There is $h\in\C[y_1,\dots,y_r][1/y_1\cdots y_r]$ with $\pi^*h=f_{\tau}$ (since $\nu^{(1)},\dots,\nu^{(r)}$ generate $V_{\tau,\Z}$), and $\rho^*y_n^{a_i}=x_i^{a_i}$. The proof of non-degeneracy for $g_{\si}$ is then reduced to that for $f_{\tau}$, using the above finite \'etale morphism $\rho$ together with Remark~(1.7), where we have the following relation of logarithmic vector fields:
$$\rho_*\bl(\msum_{j=1}^n\,c_jx_j\dd_{x_j}\br)=\msum_{k=1}^n\bl(\msum_{j=1}^n\,c_j\nu_j^{(k)}\br)y_k\dd_{y_k}.$$
This completes the proof of Proposition~(A.2).
\ms
As a corollary of Proposition~(A.2), we see that the convenience condition is not needed in the main theorem of \cite{exp}. We have more precisely the following.
\msn
{\bf Corollary~A.3.} {\it Assume $f$ has an isolated singularity with non-degenerate Newton boundary. Then the filtration on $\C\{x\}/(\dd f)$ induced by the Newton filtration $V_N^{\ssb}$ in $(1.7.4)$ coincides with the quotient filtration of the $V$-filtration on the Brieskorn lattice.}
\ms
(The argument is similar to the proof of Theorem~(2.7), using the inclusion of $(\dd f)\supset\mm_{X,0}^k$ for $k\gg 0$ together with the positivity of $\al_{f,1}$.)

\end{document}